\documentclass{article}

\usepackage[english]{babel}
\usepackage{placeins}
\usepackage[letterpaper,top=2cm,bottom=2cm,left=3cm,right=3cm,marginparwidth=1.75cm]{geometry}

\usepackage{algpseudocode}

\usepackage{algorithm}
\usepackage{breqn}
\usepackage{verbatim}
\usepackage{caption2}
\usepackage{tikz}
\usetikzlibrary{decorations.pathreplacing,decorations.markings}
\usepackage{float}
\usepackage{pifont}
\usepackage{enumerate}
\usepackage{geometry} % 调整页面布局，使表格更美观
\usepackage{amsmath}
\usepackage{graphicx}
\usepackage[colorlinks=true, allcolors=blue]{hyperref}
\newtheorem{theorem}{Theorem}[section]
\newtheorem{conjecture}[theorem]{Conjecture}
\newtheorem{definition}[theorem]{Definition}
\newtheorem{lemma}[theorem]{Lemma}
\newtheorem{proposition}[theorem]{Proposition}
\newtheorem{example}[theorem]{Example}
\newtheorem{corollary}[theorem]{Corollary}

\def\arb{arborescence}

\def\cupcup{\cup\cdots\cup}
\def\DD{\hbox{-}}
\def\CC{\hbox{-}\cdots\hbox{-}}
\def\cupcup{\cup\cdots\cup}
\def\up{\overrightarrow}

\def\qed{\hfill \rule{4pt}{7pt}}

\def\pf{\noindent {\it Proof. }} 

\title{ Optimal $\chi$-boundness of $\ell$-holed graphs}
\author{Yan Wang$^{1,}$ \thanks{Supported by National Key R\&D Program of China under Grant No. 2022YFA1006400 and Shanghai Municipal Education Commission (No. 2024AIYB003), email: yan.w@sjtu.edu.cn},  \;\; Rong Wu$^{2,}$ \thanks{Supported by National Key R\&D Program of China under Grant No. 2022YFA1006400 and Shanghai Municipal Education Commission (No. 2024AIYB003), email: rong\_w24246@163.com} \\\\
    \small $^1$School of Mathematical Sciences, CMA-Shanghai\\
    \small Shanghai Jiao Tong University, 800 Dongchuan Road, Shanghai 200240, China\\
    \small $^2$School of Mathematical Sciences \\
    \small Shanghai Jiao Tong University, 800 Dongchuan Road, Shanghai 200240, China}

\begin{document}
\maketitle

\begin{abstract}
A graph is {\em{$\ell$-holed}} if all of its induced cycles of length at least four have length exactly $\ell$. 
%A hole is said to be {\em{odd}} if it has odd length.
In the paper, we prove that if $G$ is an $\ell$-holed graph with odd $\ell\geq 7$, then $\chi(G)\leq {\lceil {\ell \over {\ell-1}}\omega(G) \rceil}$. 
This result is sharp. 
\end{abstract}

\section{Introduction } 

All graphs considered in this paper are finite, simple, and undirected. 
A graph $G$ is $k$-colorable if there exists a mapping $c: V(G)\rightarrow \{1, 2, \cdots, k\}$ such that $c(u)\neq c(v)$ whenever $uv\in E(G)$.
The \textit{chromatic number $\chi(G)$} of $G$ is the minimum integer $k$ such that $G$ is $k$-colorable.
The clique number $\omega(G)$ of $G$ is the maximum integer $k$ such that $G$ contains a complete graph of size $k$.
For a graph $G$, if $\chi(H)=\omega(H)$ for every induced subgraph $H$ of $G$, then we call $G$ a \textit{perfect} graph.
For a graph $H$, we say that $G$ is \textit{$H$-free} if $G$ has no induced subgraph isomorphic to $H$.
Let $\cal F$ be a family of graphs. We say that $G$ is \textit{$\cal F$-free} if $G$ is $F$-free for every member $F$ of $\cal F$.
If there exists a function $\phi$ such that $\chi(G)\leq \phi(\omega(G))$ for each $G\in {\cal F}$, then we say that $\cal{F}$ is $\chi$-\textit{bounded class}, and call $\phi$ a \textit{binding function} of $\cal{F}$.
The concept of $\chi$-boundedness was raised by Gy\'{a}rf\'{a}s in 1975 \cite{g2}.
Studying what families of graphs are $\chi$-bounded, and finding the optimal binding function for $\chi$-bounded class are important problems in this area.
Since the clique number is a trivial lower bound of the chromatic number, if a family of $\chi$-bounded graphs has a linear binding function, then it must be asymptotically optimal up to a constant factor.
We refer the readers to \cite{ss} for a survey on $\chi$-bounded problems.

Erd\H{o}s \cite {erdos} showed that for any positive integers $k$ and $g$, there exists a graph $G$ with $\chi (G)\geq k$ and no cycles of length less than $g$. 
This result motivates the study of the chromatic number of {$\cal{H}$}-free graphs for some {$\cal{H}$}.
Based on this, Gy\'{a}rf\'{a}s \cite{g2} and Sumner \cite{sumner} independently conjectured that if $F$ is a forest, then every $F$-free graph is $\chi$-bounded.
Due to \cite {erdos}, if no member of {$\cal{H}$} is a forest, then a necessary condition for  $\chi$-boundedness of {$\cal{H}$}-free graphs is that {$|\cal{H}|$} is infinite. 
Hence, it is natural to consider the case when $\cal{H}$ contains infinite number of induced cycles. 
%Along with the Gy\'{a}rf\'{a}s‐Sumner conjecture, Gy\'{a}rf\'{a}s\cite{g} put forward three conjectures related to holes, as follows:
%\begin{itemize}
%    \item The ideal of graphs with no odd hole is $\chi$‐bounded.
%    \item For all $\ell \geq 0$, the ideal of graphs with no hole of length $>\ell$ is $\chi$‐bounded.
%    \item For all $\ell \geq 0$, the ideal of graphs with no odd hole of length $>\ell$ is $\chi$‐bounded.
%\end{itemize}

A hole in a graph is an induced cycle of length at least $4$.
A hole is said to be \textit{odd} (resp. \textit{even}) if it has odd (resp. even) length.
Addario-Berry et al. \cite{cs08}, proved that every even hole free graph has a vertex whose neighbors are the union of two cliques, which implies that $\chi(G)\leq 2\omega(G)-1$.
However, the situation becomes much more complicated for odd hole free graphs.
The Strong Perfect Graph Theorem \cite{cs02} asserts that a graph is perfect if and only if it induces neither odd holes nor their complements.
Confirming a conjecture of Gy\'{a}rf\'{a}s \cite{g2},
Scott and Seymour \cite{sso} proved that odd hole free graphs are $\chi$-bounded with binding function $\frac{2^{2^{\omega(G)+2}}}{48(\omega(G)+2)}$.
Ho\'{a}ng and McDiarmid \cite{hm} conjectured that $\chi(G)\leq 2^{\omega(G)-1}$ for an odd hole free graph $G$.
A graph is said to be \textit{short-holed} if every hole has length 4. Sivaraman \cite{sivaraman} conjectured that $\chi(G)\leq \omega(G)^2$ for all short-holed graphs whereas the best known upper bound is $\chi(G)\leq 10^{20}2^{\omega(G)^2}$ due to Scott and Seymour \cite{ss}. 
Morever, Scott and Seymour \cite{ss2018} proved that for every integer $\ell\geq 0$, there exists $k$ such that if $G$ is triangle‐free and $\chi(G)>k$ , then $G$
has $\ell$ holes of consecutive lengths. 
In addition, Scott and Seymour \cite{ss2019} also proved that graphs containing no holes with specific residue are $\chi$-bounded.
%For all integers $k\geq 0$ and $\ell \geq 1$, the ideal of all graphs with no hole of length $k$ modulo $ell$ is $\chi$‐bounded.

A graph $G$ is $k$-$divisible$ if for every induced subgraph $H$ of $G$, either $H$ is a stable set, or the vertex set of $H$ can be partitioned into $k$ sets, none of which contains a largest clique of $H$. Hoàng and McDiarmid \cite{hm} conjectured that every odd hole free graph is $2$-divisible.
A graph $G$ is $perfectly$ $divisible$ if every induced subgraph $H$ of $G$ contains a set $X$ of vertices such that $X$ meets all largest cliques of $H$, and $X$ induces a perfect graph. Scott and Seymour \cite{ss} mentioned a conjecture of Hoàng: If a graph $G$ is odd hole free, then $V (G)$ can be partitioned into $\omega(G)$ subsets of which each induces a perfect graph.
These conjectures are known only for some special graph classes.
A $banner$ is the graph that consists of a hole on four vertices and a single vertex with precisely one neighbor on the hole.
Hoàng \cite{ct8} prove that (banner, odd hole)-free graphs are perfectly divisible. 
A $bull$ is the graph consisting of a triangle with two disjoint pendant edges.
Chudnovsky and Sivaraman \cite{ct9} proved that ($P_5,C_5$)-free graph is $2$-divisible. They also showed that either (odd hole, bull)-free or ($P_5$, bull)-free is perfectly divisible.
Chudnovsky, Robertson, Seymour and Thomas \cite{CRST10} confirmed these conjectures for $K_4$-free graphs. In fact, they showed that every $($odd hole, $K_4)$-free graph is $4$-colorable.
Recently, Sun and Wang \cite{sw} show that every $($odd hole, $\overline{2P_3})$-free graph $G$ has $\chi(G)\le \omega(G)+1$  and characterize the graphs when equality holds.

The study on the chromatic number of graphs with odd hole of exactly one length witnesses much progress recently.
The girth of a graph $G$, denoted by $g(G)$, is the minimum length of a cycle in $G$.
Let $\ell\geq 2$ be an integer. Let ${\cal G}_\ell$ denote the family of graphs that have girth $2\ell+1$ and have no odd holes of length at least $2\ell+3$. 
%The graphs in ${\cal G}_2$ are called \textit{pentagraphs}, and the graphs in ${\cal G}_3$ are called \textit{heptagraphs}.
Plummer and Zha \cite{PZ14} conjectured that every graph in ${\cal G}_2$ is $3$-colorable.
Subsequently, Chudnovsky and Seymour \cite{cs22} confirmed their conjecture.
Recently, Wu, Xu and Xu \cite{wxx22+} showed that every graph in ${\cal G}_3$ are $3$-colorable.
They also conjectured in \cite{wxx22} that all graph in $\bigcup_{\ell \ge 2}{\cal G}_\ell$ are $3$-colorable.
More recently, Chen \cite{c23} proved that all graphs in $\bigcup_{\ell\geq 5}{\cal G}_\ell$ are $3$-colorable. 
Finally, Wang and Wu \cite{ww} confirmed Wu, Xu and Xu's conjecture.

A graph is {\em{$\ell$-holed}} if all of its induced cycles of length at least four have length exactly $\ell$. 
Cook et al. \cite{chprsstv} gave a complete description of the $\ell$-holed graphs for every $\ell\geq 7$.
It is clearly that for every $\ell$-holed graph $G$ with even $\ell\geq 6$, $\chi(G)=\omega(G)$ since $G$ is perfect. 
On the other hand, $\chi(G)\leq 2\omega(G)-1$ for every $\ell$-holed graph $G$ with odd $\ell\geq 7$, since $G$ is an even hole free graph. 
However, the binding function may not be optimal for $\ell$-holed graph with odd $\ell\geq 7$.
In this paper, we prove the following. 

\begin{theorem} \label{main theorem}
For an odd integer $\ell\geq 7$, let $G$ be an $\ell$-holed graph. Then $\chi(G)\leq{\lceil {\ell \over {\ell-1}} \omega(G) \rceil}$.
\end{theorem}

Note that Theorem \ref{main theorem} improves the upper bound $2\omega(G)-1$ when $\omega(G)\geq 3$ for $\ell$-holed graphs with odd $\ell\geq 7$.
It is worthwhile to mention that Theorem \ref{main theorem} is sharp: 
consider the clique blow-up of $C_\ell$ (replacing every vertex of $C_\ell$ by a clique of arbitrary size).

Now we briefly sketch the proof of Theorem \ref{main theorem}.
By the structural description of the $\ell$-holed graphs, we know that $G$ is either a blow-up of a cycle of length $\ell$, or a blow-up of an $\ell$-framework. 
When $G$ is a blow-up of an $\ell$-cycle, we give an explicit coloring with chromatic number at most ${\lceil {\ell \over {\ell-1}} \omega(G) \rceil}$.
When $G$ is a blow-up of an $\ell$-framework, we prove by contradiction. 
Let $G$ be a minimal counterexample.
We can deduce that $m \le 4$. 
Then in each case we give a specific coloring of $G$ and show the chromatic number is at most ${\lceil {\ell \over {\ell-1}} \omega(G) \rceil}$.

The organization of this paper is as follows. 
In Section 2, we give the structure of $\ell$-holed graphs and some useful lemmas.
In Section 3, we introduce two colorings: cyclic coloring and balanced coloring,
which will be applied several times in the subsequent proof.
We present an explicit coloring to show that the blow-up of $\ell$-cycle is ${\lceil {\ell \over {\ell-1}} \omega(G) \rceil}$-colorable in Section 4.
Finally, we prove that the blow-up of $\ell$-framework is ${\lceil {\ell \over {\ell-1}} \omega(G) \rceil}$-colorable and complete the proof of Theorem \ref{main theorem} in Section 5.

\section{Preliminary}

In this section, we collect some notations and useful lemmas.
To describe the structures of $\ell$-holed graphs with odd $\ell\geq 7$ given in \cite{chprsstv}, we need the following definitions first.

An ordering of a set $X$ means a sequence enumerating the members of $X$. Let $v_1, \cdots, v_n$ be an ordering of $X\subseteq V(G)$.
We say a vertex $u\in V(G)\backslash X$ is adjacent to an initial segment of the ordering of $X$ if for all $i, j\in \{1, \cdots, n\}$ with $i<j$, if $u, v_j$ are adjacent then $u, v_i$ are adjacent. 
An ordered clique means a clique together with some ordering of it.
We will often use the same notation for an ordered clique and the (unordered) clique itself when there is no confusion on the ordering.
Let $X$ and $Y$ be disjoint subsets of $V(G)$ (with orderings). 
We denote by $G[X, Y]$ the bipartite subgraph of $G$ with vertex set $X\cup Y$ and edge set being the set of edges of $G$ between $X$ and $Y$.
%A {\em \textbf{half-graph}} is a bipartite graph with no induced two-edge matching. Take orderings $x_1, \cdots, x_m$ and $y_1, \cdots, y_n$ of $X$ and $Y$ respectively.
We say $G[X, Y]$ obeys these ordering if for all $i, i', j, j'$ with $1\leq i\leq i'\leq m$ and $1\leq j\leq j'\leq n$, if $x_{i'}y_{j'}$ is an edge then $x_iy_j$ is an edge; or equivalently, each vertex in $Y$ is adjacent to an initial segment of $x_1, \cdots, x_m$, and each vertex in $X$ is adjacent to an initial segment of $y_1, \cdots, y_n$.
%Thus, $G[X, Y]$ is a \textbf{half-graph} if and only if there are orderings of $X$ and $Y$ that $G[X, Y]$ obeys.
%If $X, Y, Z$ are disjoint ordered cliques of $G$, we say that $G[X, Y]$, $G[X, Z]$ are compatible if $G[X, Y\cup Z]$ obeys the ordering for some ordering of $Y \cup Z$.

\begin{figure}[h]
\makebox[\linewidth][c]{%
%——— 左图：9-cycle blow-up ———
\begin{minipage}{.49\linewidth}
\centering
\begin{tikzpicture}[scale=.7,auto=left]
\tikzstyle{every node}=[inner sep=1.5pt, fill=black,circle,draw]
\def\q{1.5}
\def\r{2}
\def\s{2.5}
\def\t{3}
\def\angle{360/9}
\node (a1) at ({\q*cos(90)}, {\q*sin(90)}) {};
\node (a2) at ({\r*cos(90)}, {\r*sin(90)}) {};
\node (a3) at ({\s*cos(90)}, {\s*sin(90)}) {};
\node (b1) at ({\q*cos(90+\angle)}, {\q*sin(90+\angle)}) {};
\node (b2) at ({\r*cos(90+\angle)}, {\r*sin(90+\angle)}) {};
\node (b3) at ({\s*cos(90+\angle)}, {\s*sin(90+\angle)}) {};
\node (b4) at ({\t*cos(90+\angle)}, {\t*sin(90+\angle)}) {};
\node (c1) at ({\q*cos(90+2*\angle)}, {\q*sin(90+2*\angle)}) {};
\node (c2) at ({\r*cos(90+2*\angle)}, {\r*sin(90+2*\angle)}) {};
\node (c3) at ({\s*cos(90+2*\angle)}, {\s*sin(90+2*\angle)}) {};
\node (d1) at ({\q*cos(90+3*\angle)}, {\q*sin(90+3*\angle)}) {};
\node (d2) at ({\r*cos(90+3*\angle)}, {\r*sin(90+3*\angle)}) {};
\node (e1) at ({\q*cos(90+4*\angle)}, {\q*sin(90+4*\angle)}) {};
\node (e2) at ({\r*cos(90+4*\angle)}, {\r*sin(90+4*\angle)}) {};
\node (f1) at ({\q*cos(90+5*\angle)}, {\q*sin(90+5*\angle)}) {};
\node (f2) at ({\r*cos(90+5*\angle)}, {\r*sin(90+5*\angle)}) {};
\node (f3) at ({\s*cos(90+5*\angle)}, {\s*sin(90+5*\angle)}) {};
\node (g1) at ({\q*cos(90+6*\angle)}, {\q*sin(90+6*\angle)}) {};
\node (g2) at ({\r*cos(90+6*\angle)}, {\r*sin(90+6*\angle)}) {};
\node (g3) at ({\s*cos(90+6*\angle)}, {\s*sin(90+6*\angle)}) {};
\node (h1) at ({\q*cos(90+7*\angle)}, {\q*sin(90+7*\angle)}) {};
\node (h2) at ({\r*cos(90+7*\angle)}, {\r*sin(90+7*\angle)}) {};
\node (h3) at ({\s*cos(90+7*\angle)}, {\s*sin(90+7*\angle)}) {};
\node (i1) at ({\q*cos(90+8*\angle)}, {\q*sin(90+8*\angle)}) {};
\node (i2) at ({\r*cos(90+8*\angle)}, {\r*sin(90+8*\angle)}) {};
\node (i3) at ({\s*cos(90+8*\angle)}, {\s*sin(90+8*\angle)}) {};
\node (i4) at ({\t*cos(90+8*\angle)}, {\t*sin(90+8*\angle)}) {};

\foreach \from/\to in {a1/b1,a1/b2,a1/b3,a1/b4,a2/b1,a2/b2,a2/b3,a3/b1,a3/b2}
\draw [-] (\from) -- (\to);
\foreach \from/\to in {b1/c1,b1/c2,b1/c3,b2/c1,b2/c2,b3/c1,b3/c2,b4/c1}
\draw [-] (\from) -- (\to);
\foreach \from/\to in {c1/d1,c1/d2,c2/d1,c3/d1}
\draw [-] (\from) -- (\to);
\foreach \from/\to in {d1/e1,d1/e2,d2/e1}
\draw [-] (\from) -- (\to);
\foreach \from/\to in {e1/f1,e1/f2,e1/f3,e2/f1,e2/f2}
\draw [-] (\from) -- (\to);
\foreach \from/\to in {f1/g1,f1/g2,f1/g3,f2/g1,f2/g2,f2/g3,f3/g1}
\draw [-] (\from) -- (\to);
\foreach \from/\to in {g1/h1,g1/h2,g1/h3,g2/h1,g2/h2,g3/h1}
\draw [-] (\from) -- (\to);
\foreach \from/\to in {h1/i1,h1/i2,h1/i3,h1/i4,h2/i1,h2/i2,h3/i1}
\draw [-] (\from) -- (\to);
\foreach \from/\to in {i1/a1,i1/a2,i1/a3,i2/a1,i2/a2,i2/a3,i3/a1,i3/a2,i4/a1}
\draw [-] (\from) -- (\to);

\foreach \from/\to in {a1/a2,a2/a3,b1/b2,b2/b3,b3/b4,c1/c2,c2/c3,d1/d2,e1/e2,f1/f2,f2/f3,g1/g2,g2/g3,h1/h2,h2/h3,i1/i2,i2/i3,i3/i4}
\draw [-] (\from) -- (\to);

\foreach \from/\to in {a1/a3,b1/b3,b1/b4,b2/b4,c1/c3,f1/f3,g1/g3,h1/h3,i1/i3,i1/i4,i2/i4}
\draw [bend left =30] (\from) to (\to);
\end{tikzpicture}
\caption{A blow-up of a 9-cycle.} \label{fig:cycleblowup} 
\end{minipage}
%——— 右图：7-framework ———
\begin{minipage}{.49\linewidth}
\centering
\begin{tikzpicture}[scale=.7,auto=left]
\tikzstyle{every node}=[inner sep=1.5pt, fill=black,circle,draw]

\node (a0) at (0,2) {};
\node (a1) at (1,2) {};
\node (a2) at (2,2) {};
\node (a3) at (3,2) {};
\node (a4) at (4,2) {};
\node (a5) at (6,2) {};
\node (a6) at (7,2) {};
\node (a7) at (8,2) {};
\node (a8) at (9,2) {};
\node (a9) at (10,2) {};
\node (a10) at (11,2) {};
\node (b1) at (1,0) {};
\node (b2) at (2,0) {};
\node (b3) at (3,0) {};
\node (b4) at (4,0) {};
\node (b5) at (6,0) {};
\node (b6) at (7,0) {};
\node (b7) at (8,0) {};
\node (b8) at (9,0) {};
\node (b9) at (10,0) {};
\node (b10) at (11,0) {};
\node (c1) at (1,1) {};
\node (c2) at (2,1) {};
\node (c3) at (3,1) {};
\node (c4) at (4,1) {};
\node (c5) at (6,1) {};
\node (c6) at (7,1) {};
\node (c7) at (8,1) {};
\node (c8) at (9,1) {};
\node (c9) at (10,1) {};
\node (c10) at (11,1) {};

\foreach \from/\to in {a1/c1,a2/c2,a3/c3,a4/c4,a5/c5,a6/c6,a7/c7,a8/c8,a9/c9,a10/c10}
\draw(\from) to (\to);
\foreach \from/\to in {c1/b1,c2/b2,c3/b3,c4/b4,c5/b5,c6/b6,c7/b7,c8/b8,c9/b9,c10/b10}
\draw(\from) to (\to);

\draw[fill=gray] (a0) to [bend left=20] (a9) to [bend left=10] (a10) to [bend right=23] (a0);
\draw[fill=gray] (a3) to [bend left=22] (a8) to [bend right=17] (a3);
\draw[fill=gray] (a4) to [bend left=20] (a5) to [bend left=10] (a6) to [bend left=10] (a7) to [bend right=20] (a4);

\draw[fill=gray] (b9) to [bend left=23] (b1) to [bend right=10] (b2) to [bend right=10] (b3) to [bend right=20] (b9);
\draw[fill=gray] (b8) to [bend left=20] (b4) to [bend right=10]  (b8);

\node (a0) at (0,2) {};
\node (a1) at (1,2) {};
\node (a2) at (2,2) {};
\node (a3) at (3,2) {};
\node (a4) at (4,2) {};
\node (a5) at (6,2) {};
\node (a6) at (7,2) {};
\node (a7) at (8,2) {};
\node (a8) at (9,2) {};
\node (a9) at (10,2) {};
\node (a10) at (11,2) {};
\node (b1) at (1,0) {};
\node (b2) at (2,0) {};
\node (b3) at (3,0) {};
\node (b4) at (4,0) {};
\node (b5) at (6,0) {};
\node (b6) at (7,0) {};
\node (b7) at (8,0) {};
\node (b8) at (9,0) {};
\node (b9) at (10,0) {};
\node (b10) at (11,0) {};

\tikzstyle{every node}=[]
\draw (a0) node [below right]           {\footnotesize$a_0$};
\draw (a1) node [below right]           {\footnotesize$a_1$};
\draw (a2) node [below right]           {\footnotesize$a_2$};
\draw (a3) node [below right]           {\footnotesize$a_3$};
\draw (a4) node [below right ]           {\footnotesize$a_4$};
\draw (a5) node [below right ]           {\footnotesize$a_5$};
\draw (a6) node [below right]           {\footnotesize$a_6$};
\draw (a7) node [below right]           {\footnotesize$a_7$};
\draw (a8) node [below right]           {\footnotesize$a_8$};
\draw (a9) node [below right]           {\footnotesize$a_9$};
\draw (a10) node [below right]           {\footnotesize$a_{10}$};

\draw (b1) node [above right]           {\footnotesize$b_1$};
\draw (b2) node [above right]           {\footnotesize$b_2$};
\draw (b3) node [above right]           {\footnotesize$b_3$};
\draw (b4) node [above right]           {\footnotesize$b_4$};
\draw (b5) node [above right]           {\footnotesize$b_5$};
\draw (b6) node [above right]           {\footnotesize$b_6$};
\draw (b7) node [above right]           {\footnotesize$b_7$};
\draw (b8) node [above right]           {\footnotesize$b_8$};
\draw (b9) node [above right]           {\footnotesize$b_9$};
\draw (b10) node [above right]           {\footnotesize$b_{10}$};

\begin{scope}[thick, decoration={
    markings,
    mark=at position 0.5 with {\arrow{>}}}
    ]
    \draw[postaction={decorate}] (4,3.1) to (a1);
    \draw[postaction={decorate}] (5.5,2.6) to (a4);
    \draw[postaction={decorate}] (7.3,-0.8) to (b8);
    \draw[postaction={decorate}] (6.8,-0.3) to (b7);

\foreach \from/\to in {a1/a2,a2/a3, b10/b9, b7/b6,b6/b5}
\draw[postaction={decorate}] (\from) to (\to);
\end{scope}
\end{tikzpicture}
\caption{A 7-framework with $m=4, k=10$} \label{fig:7framework} 
\end{minipage}}
\end{figure}

\begin{definition}
Let $G$ be a graph with vertex set partitioned into sets $W_1, \cdots, W_\ell$, with the following properties:
\begin{itemize}
\setlength{\itemsep}{0pt}
\setlength{\parskip}{0pt}
\setlength{\parsep}{0pt}
    \item $W_1, \cdots, W_\ell$ are non-null ordered cliques;\\
    \item for $1\leq i\leq \ell$, $G[W_{i-1}, W_i]$  obeys the ordering (reading subscripts modulo $\ell$);\\
    \item for all distinct $i, j\in \{1, \cdots, \ell\}$, if there is an edge between $W_i, W_j$ then $j \equiv i\pm 1$ (modulo $\ell$);\\
%    \item and for $1\leq i\leq \ell$, the graphs $G[W_i, W_{i+1}]$, $G[W_i, W_{i-1}]$ are compatible.
\end{itemize}
We call such a graph a {\em blow-up of an $\ell$-cycle}. (See Figure \ref{fig:cycleblowup})
\end{definition}

An {\em arborescence} is a tree with its edges directed in such a way that no two edges have a common head; or equivalently, such that
for some vertex $r(T)$ (called the {\em apex}), every edge is directed away from $r(T)$. A {\em leaf} is a vertex different from the 
apex, with outdegree zero, and
$L(T)$ denotes the set of leaves of the \arb{} $T$.

%Now we define an {\em $\ell$-framework}.
%Let us describe the key nodes and some important features of the $\ell$-framework.
\begin{definition}
For any $k\geq 3$, let $a_0,\cdots , a_{k}$ and $b_1,\cdots , b_{k}$ be vertices. 
For $1\le i\le k$, there is a path $P_i$ of length $(\ell-3)/2$ between $a_i$ and $b_i$. 
%The numbers $0\LL k$ break into two intervals $\{0\LL m\}$ and $\{m+1\LL k\}$.

Let {\em tent} be a subtree of arborescence $T$,
containing at least one leaf of $T$ (may be only one leaf), which we refer to as the \emph{bases} of a tent.
Morever, the \emph{root} of a tent is its apex.
An $\ell$-\emph{framework} consists of three main parts:
an arboresence $T$ including vertices $a_0,\cdots , a_k$, an arboresence $S$ including vertices $b_1,\cdots , b_k$, and the undirected paths $P_i$ from $a_i$ to $b_i$ 
for each $i\in [k]$.
%In the $\ell$-framework, both $T$ and $S$ include tents.
Without loss of generality, let the tents in $T$ be ``upper tents''
and those in $S$ ``lower tents''.
Let $0 \le m \le k-2$ be an integer.
The apex of each upper tent is in $\{a_0,\cdots , a_m\}$ and 
its bases are a nonempty interval of $\{a_{m+1},\cdots , a_k\}$. 
Each of $a_{m+1},\cdots , a_k$ belongs to the bases of an upper tent.
The lower tents are defined similarly.
%There can be any positive number of tents, but 
There must be a tent with apex $a_0$. 
%(There is an odd number of tents in the figure, but there could be an even number.) 
When $m=0$, there are no lower tents.
The way the upper and lower tents interleave is important; for each upper tent (except the innermost when there is an odd number of tents), 
the leftmost vertex of its base is some $a_i$, and 
$b_i$ must be the apex of some lower tent; and for each lower tent (except the innermost when there is an even number of tents), 
the rightmost vertex of its base is some $b_j$, and $a_j$ must be the apex for some upper tent. 
This gives a sort of spiral running through all the apexes of tents.

For each $i\in \{1,\cdots , m\}$, if $a_{i-1}$ is the apex of an upper tent, we call the tent $T_{i-1}$. Now, there is a directed edge
from some non-leaf vertex of $T_{i-1}$ (possibly $a_{i-1}$) to
$a_{i}$.
And if $a_{i-1}$ is not the apex of any tent, there is a directed edge from $a_{i-1}$ to $a_i$. So all these upper tents and the vertices $a_0,\cdots , a_m$, are connected up in a sequence to form one big \arb{} $T$ with apex $a_0$, and with 
set of leaves either $\{a_{m+1},\cdots , a_k\}$ or $\{a_m,\cdots , a_k\}$. 
%(In the paper, for convenience, let the set of leaves be $\{a_{m+1}\LL a_k\}$??).
%There is a directed path of $T$ that contains $a_0, a_1\LL a_m$ in order, possibly containing other vertices of $T$ between them.
Similarly for each $i\in \{m+1,\cdots , k-1\}$, if $b_{i+1}$ is the apex of a lower tent, we call the tent $S_{i+1}$.
Now, there is a directed edge
from some non-leaf vertex of $S_{i+1}$ to
$b_{i}$. And if $b_{i+1}$ is not the apex of any tent, there is a directed edge from $b_{i+1}$ to $b_i$.
So all the lower tents and the vertices $b_{m+1},\cdots , b_k$, are joined up to make one \arb{} $S$ with apex $b_k$ and with set of leaves either $\{b_1,\cdots , b_m\}$ or $\{b_1,\cdots , b_{m+1}\}$.
%(In the paper, for convenience, let the set of leaves be $\{b_{1}\LL b_{m+1}\}$??).
(See Figure \ref{fig:7framework} for a $7$-framework with $m=4$, $k=10$ where $a_0-a_9-b_9-b_3-a_3-a_8-b_8-b_4-a_4$ is a spiral running.)
\end{definition}

%Therefore, if a graph is describe as above, and each directed edge belongs to one of two \arbs{} $T,S$, each undirected edge belongs to one of the paths $P_i$. 
%We call such a graph an {\em $\ell$-framework}. (We will explain later how $\ell$-frameworks describe the structure of $\ell$-holed graphs.)
A directed path between two vertices $x, y$ means a directed path either from $x$ to $y$, or from $y$ to $x$.
Let $\ell\ge 5$ be odd and let $F$ be an $\ell$-framework, with notation as above. We observe that for $1\le i<j\le k$, either there is a directed path of $T$ between $a_i, a_j$, or there is a directed path of $S$ between $b_i, b_j$ and not both.
%\begin{theorem}\label{oddadj}

%\end{theorem}
%(To clarify: ``directed path of $T$ between $a_i, a_j$'' means a directed path either from $a_i$ to $a_j$, or from $a_j$ to $a_i$.)

\begin{definition}\label{def:l-framework}
The {\em transitive closure} $\up{T}$ of an arborescence $T$
is the undirected graph with vertex set $V(T)$ in which vertices $u,v$ are adjacent if and only if some directed path of $T$
contains both of $u,v$.
Let $F$ be an $\ell$-framework and $P_1,\cdots , P_k, T,S$ be as in the definition of an $\ell$-framework $F$.
Let
$D=\up{T}\cup \up{S}\cup P_1\cupcup P_k$. Thus
$V(D)=V(F)$, and distinct $u,v\in V(D)$ are $D$-adjacent if $u$ and $v$ are adjacent in $D$.
%either they are adjacent in some $P_i$, or there is a directed path of one of $S,T$ between $u,v$.
We say a graph $G$
is a {\em blow-up} of $F$ if
\begin{itemize}
\item $D$ is an induced subgraph of $G$, and for each $t\in V(D)$ there is a clique $W_t$ of $G$, all pairwise 
disjoint and with union $V(G)$; $W_t\cap V(D)=\{t\}$ for each $t\in V(D)$, and $W_t=\{t\}$ for each $t\in V(D)\setminus V(P_1\cupcup P_k)$.
\item For each $t\in V(D)$, there is an ordering of vertices in $W_t$ with first term $t$, say $(x_1,\cdots , x_n)$ with $x_1=t$. 
For all distinct $t,t'\in V(D)$, if $t,t'$ are not $D$-adjacent then $W_t$, $W_{t'}$ are anticomplete, and if $t,t'$ are $D$-adjacent
then $G[W_t, W_{t'}]$ obeys the ordering of $W_t,W_{t'}$, and every vertex of $G[W_t, W_{t'}]$ has positive degree.
%(Consequently, if $t,t'$ are $D$-adjacent then $t$ is complete to $W_{t'}$ and vice versa.)
\item If $t,t'\in \{a_1,\cdots , a_k\}$ or $t,t'\in \{b_1,\cdots ,b_k\}$, and $t,t'$ are $D$-adjacent, then $W_t$ is complete to $W_{t'}$.
\item For each $t\in V(T)$, if $0\le i\le m$ and $a_i,t$ are $D$-adjacent, then $W_t$ is complete to $W_{a_i}$. For each 
$t\in V(S)$, 
if  $i\in \{m+1,\cdots , k\}$ and $b_i,t$ are $D$-adjacent, then $W_t$ is complete to 
$W_{b_i}$.
\item For each upper tent $T_j$ with apex $a_j$, let $t\in L(T_j)$ and let $Q: a_0=y_1\CC y_p\DD a_j\DD z_1\CC z_q=t$ be the path of $T$ from $a_0$ to $t$. Then $W_t$ is complete to $\{y_1,\cdots , y_p, a_j\}$; $W_t$ is anticomplete to $\bigcup_{t\in T\setminus V(Q)}W_t$; and
$G[W_t, \{z_1,\cdots , z_{q-1}\}]$ obeys the ordering of $W_t$ and the ordering $z_1,\cdots , z_{q-1}$ of  $\{z_1,\cdots , z_{q-1}\}$. The same holds for lower tents with $T,a_0$ replaced by $S, b_k$.
\end{itemize}
\end{definition}

In \cite{chprsstv}, Cook et al. describe the structures of $l$-holed graphs with odd $l\geq 7$. 

\begin{lemma} \label{cook}{\em \cite{chprsstv}}
Let $G$ be a graph with no clique cutset and no universal vertex, and let $\ell \geq 7$. 
Then $G$ is an odd $\ell$-holed graph if and only if either $G$ is a blow-up of a cycle of length $\ell$, or $G$ is a blow-up of an $\ell$-framework. 
\end{lemma}

In the following part of this paper, let $G$ be a blow-up of an $\ell$-framework. We use ${\cal P}_i$ to denote the blow-up of $P_i$, and $A_i$ (respectively, $B_i$) to denote blow-up of $a_i$ (respectively, $b_i$) for $i \in \{0,\cdots , k\}$ (let $B_0=\emptyset$).
We call each ${\cal P}_i$ a \emph{clique chain}.
Let $L_{i,0}=B_i$, $L_{{\ell-1}\over 2}=A_0$ and $L_{i,j}$ be the clique in ${\cal P}_i$ at distance $j$ from $B_i$ for every $i\in [k], j\in [{{\ell-1}\over 2}-1]$.
The color set of an ordered clique is an ordered sequence of colors.
We need it to be ordered so that we can identify the color of each vertex. 
The color set of $L_{i,j}$ is denoted by $L^c_{i,j} := \{L_{i, j, 1}, L_{i, j, 2}\}$ for each $i\in [k], j\in [{{\ell-1}\over 2}-1]$.
If $j$ is even, then $|L_{i, j, 1}|=|B_i|$ and $|L_{i, j, 2}|=\omega(G)-|B_i|$.
If $j$ is odd, then $|L_{i, j, 1}|=\omega(G)-|B_i|$ and $|L_{i, j, 2}|=|B_i|$.
If $L^c_{i,j}$ has size larger than $|L_{i,j}|$, then we only use the first $|L_{i,j}|$ to color vertices in $L_{i,j}$.
Moreover, the color set of $A_0$ is denoted by $A^c_0$.

\section{Cyclic coloring and balanced coloring}

We define two colorings that will be used often in the next sections. 

\begin{definition}
    For $i \in [n]$, let $X_i$ be a vertex set. Let $Y$ be a finite set of colors. The \emph{cyclic coloring} of $\{X_1, \cdots, X_n\}$ by $Y$ is given by Algorithm \ref{alg:cyclic_coloring}.

\begin{algorithm}[h]
\caption{Ordered Cyclic Coloring with Unique Colors}
\label{alg:cyclic_coloring}
\begin{algorithmic}[1]
\Require
    \Statex Ordered family of sets $\{X_1, X_2, \dots, X_n\}$, where:
    \Statex \quad 1. $|X_1| \geq |X_2| \geq \dots \geq |X_n|$
    \Statex \quad 2. Each $X_i=\{p_{i,1}, p_{i,2}, \dots, p_{i,|X_i|}\}$  is an ordered set of distinct vertices.
    \Statex Color set $Y = \{y_1, y_2, \dots, y_{|Y|} \}$, where $|Y|\geq \sum_{i=1}^n |X_i|$.
\Ensure
    \Statex Injection $c: \bigcup_{i=1}^n X_i \rightarrow Y$.

\State \textbf{Initialization:}
\State Initialize pointers $\text{next}_i \leftarrow 1$  for all $i \in [n]$
\State Initialize global color pointer $k \leftarrow 1$ 

\State \textbf{Cyclic Coloring Process:}
\Repeat
%    \Statex \textbf{One Round:}
    \For{ $i = 1$ \textbf{to} $n$ } \Comment{Process sets in order from $X_1$ to $X_n$ }
        \If{$\text{next}_i \leq |X_i|$ } \Comment{Check if $X_i$ has uncolored points}
            \State Assign color $y_k$  to vertex $p_{i,\text{next}_i}$ (that is, $c(p_{i,\text{next}_i}) = y_k$)
            \State $\text{next}_i \leftarrow \text{next}_i + 1$  \Comment{Move to the next vertex in $X_i$ }
            \State $k \leftarrow k + 1$  \Comment{Advance to next color}
        \EndIf
    \EndFor
 %   \Statex \textbf{Termination Check:}
    \State $\text{all\_colored} \leftarrow \text{True}$ 
    \For{$i = 1 \, \textbf{ to } \, n $}
        \If{$\text{next}_i \leq |X_i|$ } 
            \State $\text{all\_colored} \leftarrow \text{False}$  
            \State \textbf{break}
        \EndIf
    \EndFor
\Until{$\text{all colored} = \text{True}$  }

\State \Return Coloring $c$
\end{algorithmic}
\end{algorithm}
\end{definition}

\begin{example}
$X_1=\{p_{1,1}, p_{1,2}, p_{1,3}\},\,X_2=\{p_{2,1}, p_{2,2}\},\, X_3=\{p_{3,1}\}$ and $Y=\{1,2,3,4,5,6,7\}$.
Then the color set of $X_1$, $X_2$ and $X_3$ are $\{1,4,6\}$, $\{2,5\}$ and $\{3\}$, respectively.
\end{example}

\begin{definition}
\label{def:balanced_coloring}
Let $H$ be an induced subgraph of a blow-up of an $\ell$-framework with $m=0$, consisting of $S$ and ${\cal P}_i$ for $i \in [k]$. 
For each $i\in [k]$, let the clique chain ${\cal P}_i:=L_{i,0}-L_{i,1}-\cdots-L_{i,n}$ with $n = \frac{\ell - 3}{2}$,  $B_i:=L_{i,0}$ be an end clique of ${\cal P}_i$ and $\omega := \omega(H)$.
Note $S=\bigcup^k_{i=1} B_i$ is a clique. We may assume that $|B_1| \ge |B_2| \ge \cdots \ge |B_k|\ge \lfloor{{\ell-1}\over 4}\rfloor\lceil{\omega\over{\ell-1}}\rceil+1$.
If $\ell\equiv 1\pmod 4$  and $|B_1|>1+{3(\ell-1)\over 8}\lceil{\omega\over{\ell-1}}\rceil$, then $k\leq 3$.
In this case, we assume that $k=3$ and 
suppose $\{B_1, B'_2, B'_3\}$ have a cyclic coloring with colors $\{1, \cdots, |B_1|+|B_2'|+|B_3'| \}$ and $\{B_2\backslash B'_2, B_3\backslash B'_3\}$ have a cyclic coloring with colors $\{ |B_1|+|B_2'|+|B_3'|+1, \cdots, |B_1|+|B_2|+|B_3| \}$, where $B'_2\subseteq B_2, B'_3\subseteq B_3$ and $|B'_2|=\lceil{{\ell-1\over 8}\lceil{\omega\over{\ell-1}}\rceil}\rceil$, $|B'_3|=\lfloor{{\ell-1\over 8}\lceil{\omega\over{\ell-1}}\rceil}\rfloor$.
Otherwise, suppose $\{B_1, \cdots, B_k\}$ have a cyclic coloring.
We may assume that %$n$ is even, and 
$|L_{i,j}|=\omega-1$ for each $i\in [k], j\in [n]$. 
Recall that $G[L_{i, j}, L_{i, j+1}]$ obeys the ordering for each $i\in [k], j\in \{0, \cdots, n-1\}$. The \emph{balanced coloring} of $H$ with $\lceil{\ell \omega\over{\ell-1}}\rceil$ colors is defined as follows: The coloring $L_{i,j}^c$ of $L_{i,j}$ consists of $L_{i,j,1}, L_{i,j,2}$ in order, which we define separately. 

\begin{enumerate}[1.]
    \item When $j\in [n]$ is even, $L_{i,j,1}$ is the colors of the first $|B_i|$ vertices in $L_{i, j}$ for each $i\in [k]$. 
    Let $0 \le t \le k-2$ be the integer such that $t \equiv {j\over 2}\lceil{\omega\over{\ell-1}}\rceil
    \pmod {k-1}$.
%    If $s_1=0$, we select the first ${{j\over 2}\lceil{\omega\over{l-1}}\rceil}\over {k-1}$ elements from each of the sets $L^c_{i+1,0},\cdots, L^c_{k,0}, L^c_{1,0}, \cdots, L^c_{i-1, 0}$ to form a color set $L^{(1)}_{i, j, 1}$.
    We select the first $\lceil{{{j\over 2}\lceil{\omega\over{\ell-1}}\rceil}\over{k-1}}\rceil$ elements from each of $L^c_{i+1,0},\cdots, L^c_{i+t,0}$, the first $\lfloor{{{j\over 2}\lceil{\omega\over{\ell-1}}\rceil}\over{k-1}}\rfloor$ elements from each of $L^c_{i+t+1,0}, \cdots, L^c_{i-1, 0}$, and first $\max\{ \lceil{{{j\over 2}\lceil{\omega\over{\ell-1}}\rceil}\over{k-1}}\rceil, |B_i| - {j\over 2}\lceil{\omega\over{\ell-1}}\rceil \}$ elements from $L^c_{i,0}$ to form a set $L^{(1)}_{i, j, 1}$.
    We sort the elements of $L^{(1)}_{i, j, 1}$ in ascending order and now $L_{i,j,1}$ are the first $|B_i|$ elements. 
    Let $L_{i,0,1}:=L^c_{i,0}$.

    When $j\in [n-1]$ is odd, $L_{i,j,1}$ is the colors of the first $\omega-|B_i|$ vertices in $L_{i, j}$ for each $i\in [k]$. 
    If $j \ge 3$ and $L_{i, j-1, 1}=\{1,\cdots,|B_i|\}$, then $L^c_{i,j}=\{\lceil{\ell\omega\over{\ell-1}}\rceil, \cdots ,\lceil{\ell\omega\over{\ell-1}}\rceil - \omega + 2 \}$.
    Otherwise, $L_{i,j,1}$ are the first $\omega-|B_i|$ elements from $\{\lceil{\ell\omega\over{\ell-1}}\rceil,\cdots,1\}\backslash (L_{i,j-1,1}\cup L_{i,j+1,1})$.
    If $j=n$ is odd, then $L_{i,j,1}$ are the first $\omega-|B_i|$ elements from $\{\lceil{\ell\omega\over{\ell-1}}\rceil,\cdots,1\}\backslash (L_{i,j-1,1}\cup \{1,\cdots,|B_i|\})$.
    
    \item Now we define $L_{i,j,2}$. 
    For each $i\in [k]$, let $J_i \in [n]$ be the smallest even integer  such that $L_{i,J_i,1}=\{1,\cdots ,|B_i|\}$ if it exists; Otherwise $J_i = n+1$.
      If $1\leq j< J_i$, when $j$ is odd, $L_{i,j,2}$ are the first $|B_i|-1$ elements in $\{\lceil{\ell\omega\over{\ell-1}}\rceil,\cdots,1\}\backslash L_{i,j,1}$;
      and when $j$ is even,
      $L_{i,j,2}$ are the first $\omega-|B_i|-1$ elements in  $\{1, \cdots,\lceil{\ell\omega\over{\ell-1}}\rceil\}\backslash L_{i,j,1}$.
     If $ J_i \le j \leq n$, when $j$ is even, 
     $L^c_{i,j}=\{1,\cdots ,\omega-1\}$;
     and when $j$ is odd,
     $L^c_{i,j}=\{\lceil{\ell\omega\over{\ell-1}}\rceil,\cdots,\lceil{\ell\omega\over{\ell-1}}\rceil - \omega + 2\}$.
     \end{enumerate}
\end{definition}

\begin{proposition}
\label{prop:balanced}
    Let $H$ be an induced subgraph of a blow-up of an $\ell$-framework. Suppose $H$ has a has a balanced coloring with $\lceil{\ell\omega\over{\ell-1}}\rceil$ colors. Then
    \begin{enumerate}[(1)]
    \item For each $i\in [n]$, even $j,j'\in\{0,\cdots,n\}$ and $j<j'$, 
    $L_{i,j,1}\backslash L^c_{i,0}\subseteq L_{i,j',1}\backslash L^c_{i,0}$. 
    \item For each $i\in [n]$, odd $j,j'\in[n]$ and $j<j'$,
    $L_{i,j,1}\cap L^c_{i,0}\subseteq L_{i,j',1}\cap L^c_{i,0}$.
    \end{enumerate}
\end{proposition}
\pf For (1), it suffices to show that $\lfloor{{{j'\over 2}\lceil{\omega\over{\ell-1}}\rceil}\over{k-1}}\rfloor \ge \lceil{{{j\over 2}\lceil{\omega\over{\ell-1}}\rceil}\over{k-1}}\rceil$ when $j'=j+2$.
For otherwise, there exists an even integer $j\in \{0,\cdots,n\}$ such that
$\lfloor{{{j+2\over 2}\lceil{\omega\over{\ell-1}}\rceil}\over{k-1}}\rfloor < \lceil{{{j\over 2}\lceil{\omega\over{\ell-1}}\rceil}\over{k-1}}\rceil$.
 Let ${j\over 2 }\lceil{\omega\over{\ell-1}}\rceil=(k-1)s_1+t_1$ and ${j+2\over 2 }\lceil{\omega\over{\ell-1}}\rceil=(k-1)s_2+t_2$ where $s_1,s_2\ge 0$ and $0 \le t_1,t_2 \le k-2$.
Obviously, $s_2\ge s_1$.
Suppose $\lfloor{{{j+2\over 2}\lceil{\omega\over{\ell-1}}\rceil}\over{k-1}}\rfloor < \lceil{{{j\over 2}\lceil{\omega\over{\ell-1}}\rceil}\over{k-1}}\rceil$.
By Definition \ref{def:balanced_coloring}, we have $t_1>t_2$, and  $t_1>0$.
Hence, $s_2 = \lfloor{{{j+2\over 2}\lceil{\omega\over{\ell-1}}\rceil}\over{k-1}}\rfloor <
\lceil{{{j\over 2}\lceil{\omega\over{\ell-1}}\rceil}\over{k-1}}\rceil=s_1+1$.
So $s_2=s_1$.
But now $t_2>t_1$, which contradicts to $t_1>t_2$.

\smallskip

For (2), suppose to the contrary that there exists an odd integer $j\in[n-2]$ such that $L_{i,j,1}\cap L^c_{i,0}\not\subseteq L_{i,j+2,1}\cap L^c_{i,0}$.
Then there exists a color $c\in (L_{i,j,1}\cap L^c_{i,0})\backslash L_{i,j+2,1}$.
Since $|L_{i,j,1}|=|L_{i,j+2,1}|$,
there exists another color $c_1\in L_{i,j+2,1}\backslash L_{i,j,1}$.
%$c,c_1\in \{1,\cdots,\lceil{\ell\omega\over\{\ell-1\}}\rceil\}$.
By Definition \ref{def:balanced_coloring},
$c\notin L_{i,j-1,1}\cup L_{i,j+1,1}$ and $c_1\notin L_{i,j+1,1}\cup L_{i,j+3,1}$ (note $L_{i,n+1,1}:=\{1,\cdots,|B_i|\}$).
By (1), either $c\notin L_{i,j+3,1}$ or $j=n-2$ and $c\in \{1,\cdots,|B_i|\}$.
By the definition of $L_{i,j+2,1}$, we have $c_1>c$.
When $c_1\in L^c_{i,0}$,
if $c_1\notin L_{i,j-1,1}$, since $c_1>c$, we have $c_1\in L_{i,j,1}$, a contradiction.
So $c_1\in L_{i,j-1,1}$.
Note $c_1>c$, then by definition $c \in L_{i,j-1,1}$, a contradiction.
When $c_1\notin L^c_{i,0}$, since $c_1\notin L_{i,j+1,1}$,
we have $c_1\notin L_{i,j-1,1}$ by (1).
Since $c_1>c$, we have $c_1\in L_{i,j,1}$, a contradiction.
\qed

\begin{example}
Let $\ell=9$ and $H$ be an induced subgraph of a blow-up of an $\ell$-framework with $m=0$ and $k=3$.
Suppose $|B_1|=17, |B_2|=12, |B_3|=11, |B'_2|=|B'_3|=5$ and $\omega=40$.
Moreover, both $\{B_1, B'_2, B'_3\}$ and $\{B_2\backslash B'_2, B_3\backslash B'_3\}$ have a cyclic coloring:  
\begin{itemize}
    \item $L^c_{1,0}=\{1,4,7,10,13,16,17,18,19,20,21,22,23,24,25,26,27\}$,\\
    $L^c_{2,0}=\{2,5,8,11,14,28,30,32,34,36,38,40\}$,\\
    $L^c_{3,0}=\{3,6,9,12,15,29,31,33,35,37,39\}$.
\end{itemize}
Then the balanced coloring of $H$ is: 
\begin{itemize}
    \item $L^c_{1,1}=\{45,\cdots,28,15,14,12,11,9,27,\cdots,16,13,10,8,7\}$,\\
    $L^c_{2,1}=\{45,\cdots,41,39,37,35,33,31,29,27,\cdots,15,13,12,10,7,40,38,36,34,32,30,28,14,11,9,8\}$,\\
    $L^c_{3,1}=\{45,\cdots,40,38,36,34,32,30,28,27,\cdots,16,14,13,11,10,8,39,37,35,33,31,29,15,12,9,7\}$.
    \item $L^c_{1,2}=\{1,\cdots,6,7,8,10,13,16,\cdots,22,9,11,12,14,15,23,\cdots,39\}$,\\
    $L^c_{2,2}=\{1,\cdots,6,8,9,11,14,28,30,7,10,12,13,15,\cdots,27,29,31,\cdots,39\}$,\\
    $L^c_{3,2}=\{1,\cdots,6,7,9,12,15,29,8,10,11,13,14,16,\cdots,28,30,\cdots,39\}$.
    \item $L^c_{1,3}=\{45,44,\cdots,7\}$,\\
    $L^c_{2,3}=\{45,\cdots,31,29,27,\cdots,16,30,28,15,14,\cdots,7\}$,\\
    $L^c_{3,3}=\{45,\cdots,30,28,\cdots,16,29,15,\cdots,7\}$.
\end{itemize}
\end{example}

\begin{lemma}\label{balanced}
Let $\ell \ge 7$ be an odd integer and $H$ be an induced subgraph of a blow-up of an $\ell$-framework with $m=0$.
A balanced coloring of $H$ with  $\lceil{\ell \omega(H) \over{\ell-1}}\rceil$ colors is a proper coloring.
\end{lemma}
\pf 
For otherwise, suppose a balanced coloring of $H$ with $\lceil{\ell \omega(H) \over{\ell-1}}\rceil$ colors is not proper.
By Definition \ref{def:balanced_coloring} we know there must exist $i\in [k]$, $j\in [n]$ such that   a color $c\in L_{i,j,2}\cap (L_{i,j-1,1}\cup  L_{i, j+1,1})$.
Note when $j=n$, $L_{i, j+1,1}$ does not exist.
%And when $j=1$, $L_{i,0,1} := L^c_{i,0}$.
Let $p_{i,j}(c)$ denote the position of color $c$ in $L^c_{i, j}$ if it exists. 
Let $L_{i,j,1}'$ be the sequence of $L_{i,j,1}$ in reverse order. 
By Definition \ref{def:balanced_coloring}, we observe that $L_{i,j,2}$ contains the subsequence of $L_{i,j-1,1}'$ before $c$,
and $L_{i,j,2}$ contains the subsequence of $L_{i,j+1,1}'$ before $c$.

%so as $p_{j+1}(c)$ and $p_{j-1}(c)$.
If $c\notin L_{i,j+1,1}\backslash L_{i,j-1,1}$,
then $c\in L_{i,j-1,1}$.
Since $L_{i,j,2}$ contains the subsequence of $L_{i,j-1,1}'$ before $c$, we have $p_{i,j}(c)\geq |L_{i,j,1}|+(|L_{i,j-1,1}|+1-p_{i,j-1}(c))=|L_{i,j,1}|+(\omega-|L_{i,j,1}|+1-p_{i,j-1}(c))$.
% So $p_{i,j}(c)+p_{i,j-1}(c) \ge \omega + 1$, a contradiction.
Since $p_{i,j+1}(c)\ge p_{i,j-1}(c)$ if $c\in L_{i,j+1,1}$,
$p_{i,j+1}(c)+p_{i,j}(c)\ge p_{i,j-1}(c)+p_{i,j}(c)\ge \omega+1$.

Otherwise, $c\in L_{i,j+1,1}\backslash L_{i,j-1,1}$,
then $0<|L_{i,j-1,1}\backslash L_{i,j+1,1}|\leq \lceil{\ell\omega\over{\ell-1}}\rceil-\omega$.
Let $b \in L_{i,j-1,1}\backslash L_{i,j+1,1}$. 
If $j$ is odd, then $b \in L_{i,0}^c$ by Proposition \ref{prop:balanced} and $b > c$ by Definition \ref{def:balanced_coloring}. So $b$ appears before $c$  in $L_{i,j,2}$ by Definition \ref{def:balanced_coloring} (2).
If $j$ is even, then $b \not\in L_{i,0}^c$ by Proposition \ref{prop:balanced}.
By Definition \ref{def:balanced_coloring}, we have $b \not\in L_{i,j,1}$, $b \in L_{i,j+2,1}$ and $c \not \in L_{i,j,1} \cup L_{i,j+2,1}$.
Thus $c > b$.
So $b$ appears before $c$ in in $L_{i,j,2}$ by Definition \ref{def:balanced_coloring} (2).
To summarize, all the colors in $L_{i,j-1,1}\backslash L_{i,j+1,1}$ appear before $c$ in $L_{i,j,2}$.
%the order of colors in $L_{i,j-1,1}\backslash L_{i,j+1,1}$ is ahead of the order of $c$ in $L_{i,j,2}$
Moreover, note that $L_{i,j,2}$ contains the sequence of $L'_{i,j+1,1}$ before $c$.
Then 
$p_{i,j}(c)\geq |L_{i,j,1}| +|L_{i,j-1,1}\backslash L_{i,j+1,1}| + (|L_{i,j+1,1}|+1-p_{i,j+1}(c))=|L_{i,j,1}| +|L_{i,j-1,1}\backslash L_{i,j+1,1}| + (\omega-|L_{i,j,1}|+1-p_{i,j+1}(c)) > \omega + 1 - p_{i,j+1}(c)$.
Hence, $p_{i,j}(c)+p_{i,j+1}(c)>\omega+1$, a contradiction.
\qed

\section{Coloring of blow-ups of a cycle of length $\ell$}

In this section, we determine the chromatic number of blow-ups of a cycle of length $\ell$.

\begin{lemma}
\label{lem:blow-up-cycle-colorable}
Let $\ell\geq 5$ be an odd integer. Every blow-up $G$ of a cycle of length $\ell$ is $\lceil \frac{\ell}{\ell - 1}w(G) \rceil $-colorable.
\end{lemma}
\pf 
Suppose $G$ is a blow-up of a cycle of length $\ell$ and $\omega = \omega(G)$. 
By definition, $V(G)$ is partitioned into sets $W_1, \cdots, W_\ell$.
For each $i\in [\ell]$, since $G[W_{i-1},W_i]$ obeys the ordering, we have $|W_i|<\omega$.
%We may assume that $|W_i| = \omega(G)-1$.
For each $i\in [\ell]$, we color each element of $W_i$ sequentially by the first $|W_i|$ colors of $X_i$, which we define explicitly in the next paragraphs. 
Let $X_i := X_{i,1}\cup X_{i,2}$ for each $i\in [\ell]$.

First suppose $\omega$ is even. 
We write ${\omega\over 2}=s{\lceil {\omega \over {\ell-1}}\rceil}+j$, where $0 \le s\leq {{\ell-3}\over 2}$ and $1\leq j\leq{\lceil {\omega \over {\ell-1}}\rceil}$ are integers.
Now we define $X_{i, 1}$.
For each $i\in [\ell]$, $X_{i, 1}$ is a sequence of colors of size ${\omega\over 2}$.
We define
$$X_{1, 1}:=\{ 1, 2, \cdots, {\omega\over 2} \},$$ 
$$X_{2, 1}:=\{ {\omega\over 2}+1,  {\omega\over 2}+2, \cdots, \omega \},$$ 
$$X_{3,1}:=\{ {\lceil {\omega \over {\ell-1}}\rceil}+1,\cdots, {\omega\over 2}, \omega+1, \cdots, \omega+{\lceil {\omega \over {\ell-1}}\rceil} \},$$
and for each $2 \le h\leq s$, 
$$X_{2h, 1}:=\{ {\omega\over 2}+1+(h-1){\lceil {\omega \over {\ell-1}}\rceil}, \cdots, \omega, 1, \cdots, (h-1){\lceil {\omega \over {\ell-1}}\rceil} \},$$
%and for each $2 < h\leq s+1$, 
$$X_{2h+1, 1}:=\{ {h{\lceil {\omega \over {\ell-1}}\rceil}+1,\cdots, {\omega\over 2}, \omega+1, \cdots, \omega+{\lceil {\omega \over {\ell-1}}\rceil}, {\omega\over 2} +1},{\cdots, {\omega\over 2}+(h-1){\lceil {\omega \over {\ell-1}}\rceil}} \},$$
and
$$X_{2s+2, 1}:=\{\omega+1-j, \cdots, \omega, 1, \cdots, {\omega \over 2}-j\},$$
$$X_{2s+3, 1} :=\{ \omega+{\lceil {\omega \over {\ell-1}}\rceil}-j+1, \cdots, \omega+{\lceil {\omega \over {\ell-1}}\rceil}, {\omega \over 2}+1, \cdots, \omega -j\}.$$ 
For each $2s+2\leq m\leq \ell$, $X_{m, 1}:=X_{2s+2, 1}$ for even $m$ and $X_{m, 1}:=X_{2s+3, 1}$ for odd $m$.

Next we define $X_{i, 2}$.
Let $X_{i, 2}:=X'_{i+1, 1}$ for each $i\in [\ell]$  where $X'_{i+1, 1}$ is the sequence of $X_{i+1, 1}$ in reverse order. 
%and $X'_{\ell+1, 1}=X'_{1, 1}$.
%Let $X_i=X_{i, 1}\cup X_{i, 2}$ for each $i\in [\ell]$.
Now we show such coloring is proper.
Let $c$ be a color in $X_{i,2}$ and $p_t(c)$ denote the position of color $c$ in $X_{t}$ for $t \in [\ell]$ if it exists.
So $p_i(c)=\omega+1-p_{i+1}(c)$.
We may assume that $c\in X_{i-1,1}$. 
By construction, we have $p_{i-1}(c)\ge p_{i+1}(c)$.
%If $c\notin X_{i-1,1}$, $p_{i-1}(c)$ does not exist.
Since $p_{i-1}(c)+p_{i}(c)\ge p_{i+1}(c)+p_i(c)=\omega+1$,
such coloring is proper.
Therefore, we have $\chi(G)\leq {\lceil {\ell \omega \over {\ell-1}}\rceil}$.

Now suppose $\omega$ is odd.
Let $G'$ be the graph obtained from $G$ by deleting all the vertices from position ${\omega+1} \over 2$ to position $|W_i|$ in each $W_i$ for $i\in [\ell]$.
Note that $\omega(G') \le \omega - 1$.
By the above coloring, it is easy to see that $\chi(G') \le \lceil {\ell \over {\ell-1} }\omega(G')\rceil \le \lceil {\ell  \over {\ell-1}}(\omega-1)\rceil$.
Now in $G$, we color the vertex in position ${\omega+1} \over 2$ in each $W_i$ for $i\in [\ell]$ by the same new color, and color the vertices after position ${\omega+1} \over 2$ in $W_i$ by the reverse order of the colors of the first $\omega-1\over 2$ vertices in $W_{i+1}$ in $G'$.
Since $\lceil {\ell (\omega -1)\over {\ell-1}}\rceil+1 \le \lceil {\ell \omega \over {\ell-1}}\rceil$, we have $\chi(G)\leq {\lceil {\ell \omega \over {\ell-1}}\rceil}$.\qed

\section{Coloring of blow-up of $\ell$-frameworks}

In this section, we prove the following.

\begin{lemma}
\label{lem-color-framework}
Every blow-up of $\ell$-frameworks $G$ is $\lceil \frac{\ell}{\ell - 1}w(G) \rceil $-colorable.
\end{lemma}

By Definition \ref{def:l-framework}, 
let $A=\cup_{t\in \{a_0,\dots,a_k\}}W_t$ and $B=\cup_{s \in \{b_1,\dots,b_k\} } W_s,$
let $A^{(1)} = \cup_{t \in \up{T} \setminus \{a_0,\dots,a_k\} } W_t$
and $B^{(1)} = \cup_{s \in \up{S} \setminus \{b_1,\dots,b_k\} } W_s$.
%in addition to $W_{a_i}$, for each $i\in \{0,\cdots,k\}$, there are the other vertices in the blow up of $\up{T}$, so as the $\up{S}$. We denote these two sets of vertices as $A'$ and $B'$, respectively.
For each $i\in[k]$, we also denote $W_{a_i}$ and $W_{b_i}$ by $A_i$ and $B_i$, respectively.
Clearly, each vertex in $A^{(1)}$ (respectively $B^{(1)}$) must have a neighbor in some $A_i$ (respectively $B_i$) with $i\in [k]$.
By Definition \ref{def:l-framework}, $A^{(1)}$ is complete to $\{a_0\}$.
If a vertex $v \in A^{(1)}$ has a neighbor in $\bigcup^m_{i=1}A_i$, then $v$ is complete to $\bigcup^m_{i=1}A_i$.
Moreover, every vertex in $A^{(1)}$ has a neighbor in $\bigcup^k_{i=m+1}A_i$. 
So when $G$ is a minimum counterexample of Lemma \ref{lem-color-framework}, for $v\in A^{(1)}$,  there exists two distinct $j, j'\in [k]$ such that $A_j$ is anticomplete to $A_{j'}$ and $N_{A_j}(v)\neq \emptyset$, $N_{A_{j'}}(v)\neq \emptyset$.
%$A_{0,h}:=A_0\cup A'_h$ for each $h\in[k-m]$
%{\color{red}Let $A_0:=\{a_0\}\cup A'$ and
%B'_k=B_k\cup B'$.
%Now $A_0$ and $A_i$ obey the orderings for $i \in [k]$, and similarly $B'_k$ and $B_i$ obey the orderings for $i \in [k-1]$.}
\begin{comment}
{\color{red} Let $A''\subseteq A'$ denote the set of common interior vertices on the path from $a_0$ to  $a_i$ for each $i\in[k]$.
And let $A_0:=\{a_0\}\cup A''$, 
then $A_0$ is complete to each $A_i$.
\end{comment}

%NOT NEEDED
%Suppose $G$ is a minimum counterexample of Lemma \ref{lem-color-framework} with minimum $|V(G)|$.
%Then each vertex in $A'$ must have neighbors in at least two different ${A_i}$'s with $i \in [k]$. (EXPLAIN) 
\begin{comment}

{\color {blue}For $i \in \{m+1,\dots,k\}$, let $A'_h$ be a maximal ordered clique such that $v\in A'$ has a neighbor in $A_h$ and has no neighbors in $\bigcup^{k}_{i=h+1}A_i$.
For $i \in \{1,\dots,m\}$, let $B'_h$ be a maximal ordered clique such that $v\in B'$ has a neighbor in $B_h$ and has no neighbors in $\bigcup^{h-1}_{i=1}B_i$.
So $A'=\bigcup^{k}_{h=m+1} A'_h$ and $B'=\bigcup^{m}_{h=1} B'_h$. 
Note that it is possible that for distinct $h,h' \in \{m+1,\cdots,k\}$, $A'_h\cap A'_{h'}\neq \emptyset$, and for distinct $h,h' \in \{1,\cdots,m\}$, $B'_h\cap B'_{h'}\neq \emptyset$.}
\end{comment}

Suppose $G$ is a minimum counterexample of Lemma \ref{lem-color-framework} with minimum $|V(G)|$.
Let $s$ be a positive integer such that either $\ell=4s+3$ and $s \ge 1$ or $\ell=4s+1$ and  $s\geq 2$.
We say an ordered set $Z=\{z_1,\cdots,z_k\}$  has a coloring $c$ if $z_i$ is colored by $c(z_i)$  for every $i\in[k]$, and its color set is $c(Z):= \{c(z_1),\cdots,c(z_k)\}$. 
If we say the color set of an ordered set $Z$ is $C$, then this means we only know the color set without knowing the colors of each element.
First we show the following.

%Next, we show that if $G$ is a blow-up  of an $\ell$-framework, either $\ell=4s+3$ and  $s$ is a positive integer or $\ell=4s+1$ and  $s\geq 2$ is a integer, then $\chi(G)\leq {\lceil {\ell \omega(G)\over {\ell-1}}\rceil}$. 
%Suppose it is not true. Let $G$ be a minimal counterexample.
%It is clearly that the blow-up of $\up{T}$ are consisted of $A_i$ for each $i\in \{0,\cdots ,k\}$ and the blow-up of $\up{S}$ are consisted of $B_j$ for each $j\in \{1,\cdots , k\}$.

\begin{lemma} \label{lem11}
Let $m\geq 1$ and $P := Z_0Z_1\cdots Z_m$ be the blow-up of a path where $Z_i$ is an ordered clique for each $i\in \{0,\cdots, m\}$ and
$G[Z_j, Z_{j+1}]$  obeys the ordering for each $j\in \{0,\cdots, m-1\}$.
Let $\omega := \omega(P)$ and $\chi'>\omega$ be an integer.
Suppose $|Z_0|\leq {\omega\over 2}$.
Let $C_0$ and $C_m$ denote the sets of colors  of $Z_0$ and $Z_m$, respectively.
Moreover, suppose that at most one of $Z_0$ and $Z_m$ has a coloring $c$.
\begin{enumerate}[(1)]
\item If $m$ is odd and $|C_0\cap C_m|\leq {(m-1)(\chi'-\omega)\over 2}$, then $P$ is $\chi'$-colorable. 
\item  If $m$ is even and $|C_0\backslash C_m|\leq {m(\chi'-\omega)\over 2}$ or $|C_m\backslash C_0|\leq {m(\chi'-\omega)\over 2}$, then $P$ is $\chi'$-colorable.  
 \end{enumerate}
\end{lemma}
\pf 
We construct a proper coloring of $\chi'$ colors that satisfies (1) (respectively, (2)).
For each $i\in [m]$, since $G[Z_{i-1},Z_i]$ obeys the ordering, we have $|Z_i|<\omega$.
%Assume that except $Z_0$ and $Z_m$, for each $j\in [m-1]$, $|Z_j|=\omega-1$.
%We may assume that $|W_i| = \omega(G)-1$.
Let $X_i := X_{i,1}\cup X_{i,2}$ for each $i\in [m-1]$ and
let $|X_{i,1}|=|Z_0|, |X_{i,2}|=\omega-|Z_0|$  for every even $i\in [m-1]$ and
$|X_{i,1}|=\omega-|Z_0|, |X_{i,2}|=|Z_0|$  for every odd $i\in [m-1]$.
For each $i\in [m-1]$, we color each element of $Z_i$ sequentially by the first $|Z_i|$ colors of $X_i$.
Without loss of generality, we may assume that $|Z_0|\le |Z_m|$.

For (1), we may assume that $m \ge 3$
and $x:=|C_0\cap C_m|$.
If $Z_0$ has a coloring $c$,
by permuting the colors,
we may assume $c(Z_m)=\{\chi',\cdots, \chi'+1-|Z_m|\}$
and $C_0 = Y_1 \cup Y_2$ where $Y_1:=\{1,\cdots, |Z_0|-x\}$ and $Y_2:=\{\chi',\cdots, \chi'+1-x\}$.
% $Y_1$ and $Y_2$ is ordered set
That is, $c(Z_0)=\{c_1,c_2,\cdots,c_{|Z_0|}\}$ such that for distinct $1 \le i<j\le |Z_0|$, if $c_i, c_j\in Y_1$ then $c_i<c_j$,
and if $c_i, c_j\in Y_2$ then $c_i>c_j$.
Then we define $X_{m-1}=\{1,\cdots,\omega\}$ and $X_{m-1,1}=\{1,\cdots,|Z_0|\}$.
For each $1\leq s_1\le{m-1\over 2}$,
if $s_1(\chi'-\omega)< x$, 
then let $Y=\{c_1,\cdots,c_y\}\subseteq c(Z_0)$ where $c_y=\chi'+1-s_1(\chi'-\omega)$ and
let $Y'$ denote the ordered set consisting of the first $|Z_0|-|Y|$ elements of $\{1,\cdots,|Z_0|\}\backslash Y$ and
$X_{m-1-2s_1,1}=Y\cup Y'$;
otherwise, $X_{m-1-2s_1,1}=c(Z_0)$.
Let $X_{0,1}:=c(Z_0)$ and $X_{m-2s_1,1}$ be the ordered set consisting of the first $\omega-|Z_0|$ elements of $\{\chi',\cdots,1\}\backslash (X_{m-2s_1+1,1}\cup X_{m-2s_1-1,1})$.
And for each $i\in [m-2]$, $X_{i,2}:=X'_{i+1,1}$ where $X'_{i+1,1}$ is the sequence of $X_{i+1,1}$ in reverse order.
One can verify that this coloring is proper
and $P$ is $\chi'$-colorable.

If $Z_0$ does not have a coloring $c$,
by permuting the colors,
we may assume $c(Z_m)=\{\chi',\cdots, \chi'+1-|Z_m|\}$
and $c(Z_0)=\{1,\cdots,|Z_0|-x\}\cup \{c_1, c_2,\cdots, c_x\}$ where $c_1<c_2\cdots<c_x$ and $\{1,\cdots,|Z_0|-x\}\cap \{c_1, c_2,\cdots, c_x\}=\emptyset$.
Then we define $X_{m-1}=\{1,\cdots,\omega\}$ and $X_{m-1,1}=\{1,\cdots,|Z_0|\}$.
For each $1\leq s_1\le{m-1\over 2}$,
if $s_1(\chi'-\omega)< x$, 
then let $Y=\{1,\cdots,|Z_0|-x+s_1(\chi'-\omega)\}$ and
let $Y'$ denote the ordered set consisting of the first $|Z_0|-|Y|$ elements of $\{c_1, c_2,\cdots, c_x\}\backslash Y$ and
$X_{2s_1,1}=Y\cup Y'$;
otherwise, $X_{2s_1,1}=X_{m-1,1}$.
Let $X_{0,1}:=c(Z_0)$ and $X_{2s_1-1,1}$ denote the ordered set consisting of the first $\omega-|Z_0|$ elements of $\{\chi',\cdots,1\}\backslash (X_{2s_1,1}\cup X_{2s_1-2,1})$.
And for each $i\in [m-2]$, $X_{i,2}:=X'_{i-1,1}$ where $X'_{i-1,1}$ is the sequence of $X_{i-1,1}$ in reverse order.
One can verify that this coloring is proper
and $P$ is $\chi'$-colorable.

\begin{comment}
For (2), 
when $m\geq 4$, by permuting the colors, let $c(Z_m)=\{1,\cdots,|Z_m|\}$ and $X_{m-1}=\{\chi',\cdots,\chi'+2-\omega\}$, $X_{m-1,1}=\{\chi',\cdots,\chi'+1-(\omega-|Z_0|)\}$.
By (1), One can verify that this coloring is proper
and $P$ is $\chi'$-colorable.

When $m=2$, by permuting the colors, let $c(Z_m)=\{1,\cdots,|Z_m|\}$, $X_{m,1}=\{1,\cdots,|Z_0|\}$
and $x:=|C_0\backslash C_m|$.
If $Z_0$ has a coloring $c$
we may assume  $C_0 = Y_1 \cup Y_2$ where $Y_1:=\{1,\cdots, |Z_0|-x\}$ and $Y_2:=\{\chi',\cdots, \chi'+1-x\}$.
% $Y_1$ and $Y_2$ is ordered set
That is, $X_{0,1}:=c(Z_0)=\{c_1,c_2,\cdots,c_{|Z_0|}\}$ such that for distinct $1 \le i<j\le |Z_0|$, if $c_i, c_j\in Y_1$ then $c_i<c_j$,
and if $c_i, c_j\in Y_2$ then $c_i>c_j$.
Let $X_{1,1}$ be the ordered set consisting of the first $\omega-|Z_0|$ elements of $\{\chi',\cdots,1\}\backslash (X_{0,1}\cup X_{2,1})$.
And $X_{1,2}:=X'_{2,1}$ where $X'_{2,1}$ is the sequence of $X_{2,1}$ in reverse order.
One can verify that this coloring is proper
and $P$ is $\chi'$-colorable.
If $Z_0$ does not have a coloring $c$,
we may assume  $X_{0,1}:=c(Z_0)=\{1,\cdots,|Z_0|-x\}\cup \{c_1, c_2,\cdots, c_x\}$ where $c_1<c_2\cdots<c_x$ and $\{1,\cdots,|Z_0|-x\}\cap \{c_1, c_2,\cdots, c_x\}=\emptyset$.
Let $X_{1,1}$ be the ordered set consisting of the first $\omega-|Z_0|$ elements of $\{\chi',\cdots,1\}\backslash (X_{0,1}\cup X_{2,1})$.
And $X_{1,2}:=X'_{0,1}$ where $X'_{0,1}$ is the sequence of $X_{0,1}$ in reverse order.
One can verify that this coloring is proper
and $P$ is $\chi'$-colorable.
\end{comment}
For (2), we construct a path $P'$ such that $P':=Z_0Z_1\cdots Z_mZ_{m+1}$ and $Z_{m+1}$ is an ordered clique and complete to $Z_m$.
By permuting the colors, let $c(Z_m)=\{1,\cdots,|Z_m|\}$ and $c(Z_{m+1})=\{\chi',\cdots,\chi'+1-|Z_{m+1}|\}$.
Since $|C_0\cap C_{m+1}|\leq |C_0\backslash C_m|\leq {m(\chi'-\omega)\over 2}$,
when $Z_0$ has a coloring $c$,
by (1), $P'$ is $\chi'$-colorable, then $P$ is $\chi'$-colorable.
When $Z_0$ dose not have a coloring $c$, 
let $Z_m$ be uncolored,
by (1), $X_{m}=\{1,\cdots,|Z_0|\}\cup X'_{m-1,1}$.
But now, $c(Z_m)=\{1,\cdots, |Z_m|\}$, 
one can verify that this coloring is proper
and $P$ is $\chi'$-colorable.
\qed

\begin{lemma}\label{Claim 3.1.1}
Let $G$ be a minimal counterexample of Lemma \ref{lem-color-framework}.
%Let $A_{0,h}:=\{a_0\}\cup A'_h$ for each $h\in\{m+1,\cdots,k\}$ and $B_{k,h}:=B_k\cup B'_h$ for each $h\in\{1,\cdots,m\}$.
%Obviously, $A_{0,h}$ is a subset of $A_0$ and $B_{k,h}$ is a subset of $B_k\cup B'$.
Then we have the following.
    \begin{enumerate}[(1)]
         \item When $m=0$, $|B_j|\geq \min \{ \lceil {\omega\over {l-1}}\rceil s, \frac{\omega}{2} \}+1$ for each $j\in [k]$
         %and $|A_{0,h}|\geq \min \{ \lceil {\omega\over {l-1}}\rceil s, \omega/2 \}+1$.
         \item When $m\neq0$, $|A_i|\geq \min \{ \lceil {\omega\over {l-1}}\rceil s, \frac{\omega}{2} \}+1$ for each $i\in [m+1]$;
         $|B_j|\geq \min \{ \lceil {\omega\over {l-1}}\rceil s, \frac{\omega}{2} \}+1$ for each $j \in \{m+1,\cdots , k\}$
       %{\color{red}and if $A_{0,h}$ is complete to $\bigcup^m_{i=1}A_i$,
        % then
       %  $|A_{0,h}|\geq \min \{ \lceil {\omega\over {l-1}}\rceil s, \frac{\omega}{2} \}+1$.}   
    \end{enumerate}
\end{lemma}
\pf When $m=0$, we know that $B^{(1)}=\emptyset$ and
$B_i$ is complete to $B_j$ for every distinct $i,j\in[k]$ by Definition \ref{def:l-framework}.
Suppose there exists $j\in [k]$ such that $|B_j|\leq \min\{\lceil {\omega \over {\ell-1}}\rceil s,\ \frac{\omega}{2}\}$.
We assume that $|B_1|\leq \min\{\lceil {\omega \over {\ell-1}}\rceil s,\ \frac{\omega}{2}\}$.
And there exists a maximum subset $A'$ of $A^{(1)}$ such that $G[\{a_0\}\cup A', A_1]$ obeys the ordering.
Then $G\backslash({\cal P}_1 \setminus  B_1)$ is $\lceil {\ell\omega \over {\ell-1}}\rceil$-colorable since $G$ is a minimal counterexample.
Since $B_1$ is complete to $B_j$ for $j\in \{2,\cdots, k\}$, 
$B_1$ has a color set, but not a specific coloring as we can permute the colors of vertices in $B_1$.
By Lemma \ref{lem11} with $(Z_0,\dots,Z_m)_{\ref{lem11}} = (B_1,\cdots, \{a_0\}\cup A')$, $G$ is $\lceil {\ell\omega \over {\ell-1}}\rceil$-colorable, a contradiction. This proves (1).

When $m\neq 0$, suppose there exists $i\in [m+1]$ such that $|A_i|\leq \min\{\lceil {\omega \over {\ell-1}}\rceil s,\ \frac{\omega}{2}\}$.
Then $G\backslash({\cal P}_i\backslash A_i)$ is $\lceil {\ell\omega \over {\ell-1}}\rceil$-colorable since $G$ is a minimum counterexample.
When $i\in [m]$, by Definition \ref{def:l-framework}, $A_i$ has a color set but not a specific coloring.
When $i=m+1$, by Definition \ref{def:l-framework},
$N_{B\cup B^{(1)}}(B_i)$ is complete to $B_i$, so $N_{B\cup B^{(1)}}(B_i)$ does not have a specific coloring as we can permute the colors of vertices in $N_{B\cup B^{(1)}}(B_i)$.
Then by Lemma \ref{lem11} with $(Z_0,\dots,Z_m)_{\ref{lem11}} = (A_i,\cdots, N_{B\cup B^{(1)}}(B_i))$, $G$ is $\lceil {\ell\omega \over {\ell-1}}\rceil$-colorable, a contradiction.
Similarly, we have $|B_j|\geq \min \{ \lceil {\omega\over {l-1}}\rceil s, \frac{\omega}{2} \}+1$ for each $j\in \{m+1,\cdots , k\}$.
\begin{comment}
Now assume that $A_{0,h}$ is complete to $\bigcup^m_{i=1}A_i$.
Suppose for contradiction that $|A_{0,h}|\leq \min\{\lceil {\omega \over {\ell-1}}\rceil s,\ \frac{\omega}{2}\}$. 
Then $G\backslash({\cal P}_k\backslash B_k)$ is $\lceil {\ell\omega \over {\ell-1}}\rceil$-colorable since $G$ is a minimum counterexample. 
By Lemma \ref{lem1} with $(Z_0,\dots,Z_m)_{\ref{lem1}} = (A_{0,h},\cdots, B'_{k})$, $G$ is $\lceil {\ell\omega \over {\ell-1}}\rceil$-colorable, a contradiction.
\end{comment}
\qed

\subsection{The case when $\ell \equiv 3 \pmod 4$}

Since both $\bigcup^{m+1}_{i=1}A_i$ and $\bigcup^{k}_{j=m+1}B_j$ are cliques,
if $\ell=4s+3$ and  $s$ is a positive integer, by Lemma \ref{Claim 3.1.1},
either $\ell=7$, $m$ is at most $4$, $k-m$ is at most $5$; or $\ell\geq 11$, $m$ is at most $3$ and $k-m$ is at most $4$.
In fact, in terms of structure, $S$ and $T$ are symmetrical.
So we assume that $k-m\geq m$. 
%since when $k-m=m+1$, the structure of the $l$-holed graphs is the same as the structure of the case with $m=0$.

\begin{lemma}
    $m\neq0.$
\end{lemma}

\pf For otherwise, suppose $m = 0$. 
Let $m_i:=|B_i|$ for each $i\in [k]$.
We may assume that $m_1\geq m_2\geq \cdots \geq m_k\geq {\lceil {\omega \over {\ell-1}}\rceil}s+1$.
And we assume that $|A_0|=\omega-1$,
let $A_0$ and $A_i$ obey the ordering for each $i\in [k]$.
When $l=7$, $3\leq k\leq 5$; and when $l\geq 11$, $3\leq k\leq 4$.
Let $i'\in [k]$ be maximum such that $m_{i'}\geq \lceil{sk\lceil{\omega\over{\ell-1}}\rceil\over {k-1}}\rceil$.
If $i'$ does not exist, then $m_1< \lceil{sk\lceil{\omega\over{\ell-1}}\rceil\over {k-1}}\rceil$, and we define $i'=0$.

If $i'=0$, let $\{B_1,\cdots, B_k\}$ have a cyclic coloring and let $G\backslash A_0$ have a balanced coloring.
If $i'\neq 0$, let $\{B_1,\cdots, B_{i'}, B'_{i'+1}, \cdots, B'_k\}$ have a cyclic coloring and let $\{B_{i'+1}\backslash B'_{i'+1}, \cdots, B_k\backslash B'_k\}$ have a cyclic coloring, %where $q$ is largest number less than $i'$ such that $|B_{q}|>\lceil{sk\lceil{\omega\over{l-1}}\rceil\over {k-1}}\rceil$ and for each $i\in \{q+1,\cdots,k\}$, $|B'_i|=\lceil{s\lceil{\omega\over{l-1}}\rceil\over {k-1}}\rceil$.
and $G\backslash A_0$ have a balanced  coloring.

Now, if $i'=0$, let $A_0=\{\lceil{\ell\omega\over {\ell-1}}\rceil,\cdots, \lceil{\ell\omega\over {\ell-1}}\rceil-\omega+2\}$.
Obviously, $G$ is $\lceil{\ell\omega\over {\ell-1}}\rceil$-colorable.
If $i'\neq 0$, let $A_0:=A_{0,1}\cup A_{0,2}$
where $|A_{0,1}|=\omega-m_1$.
When $m_1\geq |A_i|>m_i$ for each $i\in [k]$, the colors of the the $(m_i+1)$-th vertex to the $m_1$-th vertex in $A_i$ already exist in other $L_{j,2s,1}$, $j\in [k]\backslash\{i\}$.
When $|A_i|>m_1$ for each $i\in [k]$, we can use at most one new color to color $A_i$, which does not contradict $A^c_{0,1}$.
Now we show that for each $3\leq k\leq 5$, if $|A_i|=m_i$ for each $i\in [k]$, then $G\backslash A_{0,2}$ is $\lceil{\ell\omega\over {\ell-1}}\rceil$-colorable. 
So it suffices to show 
\begin{equation}\label{3.1}
\sum^{i'}_{j=1}(m_j-\lceil{sk\lceil{\omega\over{\ell-1}}\rceil\over {k-1}}\rceil)+k\lceil{s\lceil{\omega\over{\ell-1}}\rceil\over {k-1}}\rceil+\omega-m_1\leq \lceil{\ell\omega\over{\ell-1}}\rceil.
\end{equation}
The proof of (\ref{3.1}) 
involves computation by cases, which we postpone to the Appendix. 
Let $A^c_{0,2}=\{\lceil{\ell\omega\over{\ell-1}}\rceil,\cdots,1\}\backslash A^c_{0,1}$ and $A^c_{0}=A^c_{0,1}\cup A^c_{0,2}$.
It is easy to see that $G$ is $\lceil{\ell\omega\over{\ell-1}}\rceil$-colorable.\qed
\begin{comment}
\begin{figure}[htp]
    \centering
    \includegraphics[width=0.38\linewidth]{5.5.png}
    \caption{$G'$ with $\ell=7, k=3$}
    \label{f2,3}
\end{figure}
\end{comment}
\begin{figure}[htbp]
\centering
% ---------- 左图 ----------
\begin{minipage}[t]{0.48\linewidth}
\centering
\begin{tikzpicture}[
  main_node/.style={circle,draw,minimum size=7mm,inner sep=2pt},
  big_node/.style={circle,draw,minimum size=10mm,inner sep=3pt}
]

% ---------------- 坐标定义 ----------------
\def\xA{-2.75}  \def\yA{4.85}   % A0
\def\xB{-4.85}  \def\yB{3.65}   % A1
\def\xC{-2.75}  \def\yC{3.65}   % A2
\def\xD{-0.60}  \def\yD{3.65}   % A3

\def\xE{-4.85}  \def\yE{2.35}   % 4
\def\xF{-2.75}  \def\yF{2.35}   % 5
\def\xG{-0.60}  \def\yG{2.35}   % 6

\def\xH{-4.85}  \def\yH{1.00}   % B1
\def\xI{-2.75}  \def\yI{1.00}   % B2
\def\xJ{-0.60}  \def\yJ{1.00}   % B3

% ---------------- 节点 ----------------
\node[big_node] (0) at (\xA,5.4)  {$A_0$};
\node[big_node] (1) at (\xB,3.6)  {$A_1$};
\node[big_node] (2) at (\xC,3.6)  {$A_2$};
\node[big_node] (3) at (\xD,3.6)  {$A_3$};

\node[big_node] (4) at (\xE,1.8)  {};
\node[big_node] (5) at (\xF,1.8)  {};
\node[big_node] (6) at (\xG,1.8)  {};

\node[big_node] (7) at (\xH,0)  {$B_1$};
\node[big_node] (8) at (\xI,0)  {$B_2$};
\node[big_node] (9) at (\xJ,0)  {$B_3$};

\node[main_node]  (10) at (-1.6,5.4) {$A'$};

% ---------------- 边 ----------------
% 红色边
\path[draw, red, thick]
(10) edge (1)
(10) edge (2)
(10) edge (3)
(1)  edge (4)
(2)  edge (5)
(3)  edge (6)
(4)  edge (7)
(5)  edge (8)
(6)  edge (9);

% 黑色边
\path[draw, black, thick]
(0) edge (1) 
(0) edge (2) 
(0) edge (3) 
(7) edge (8) 
(8) edge (9) 
(7) edge[bend right=35] (9)  % 曲线
(0) edge (10);
\end{tikzpicture}
\caption{$G'$ with $\ell=7$ and $k=3$}
\label{f2,3}
\end{minipage}
\hfill
% ---------- 右图 ----------
\begin{minipage}[t]{0.48\linewidth}
\centering
\begin{tikzpicture}[
    main_node/.style={circle, draw, minimum size=10mm, inner sep=3pt},
    small_node/.style={circle, draw, minimum size=7mm, inner sep=2pt},
    red_edge/.style={draw, thick, red}
]

% 定义水平位置
\def\xA{-4.5} % 1,5,9,13
\def\xB{-2.5} % 2,6,10,14
\def\xC{-0.5} % 3,7,11,15
\def\xD{1.5}  % 4,8,12,16

% 定义垂直位置
\def\yA{3.9}  % 0,1,2,3,4
\def\yB{2.6}  % 5,6,7,8
\def\yC{1.3}  % 9,10,11,12
\def\yD{0.0}  % 13,14,15,16
%\yE 不再使用

% 绘制节点
\node[main_node] (0) at (-0.5, 5.2) {$A_0$}; % 点0位置略高
\node[main_node] (1) at (\xA, \yA) {$A_1$};
\node[main_node] (2) at (\xB, \yA) {$A_2$};
\node[main_node] (3) at (\xC, \yA) {$A_3$};
\node[main_node] (4) at (\xD, \yA) {$A_4$};
\node[main_node] (5) at (\xA, \yB) {};
\node[main_node] (6) at (\xB, \yB) {};
\node[main_node] (7) at (\xC, \yB) {};
\node[main_node] (8) at (\xD, \yB) {};
\node[main_node] (9) at (\xA, \yC) {};
\node[main_node] (10) at (\xB, \yC) {};
\node[main_node] (11) at (\xC, \yC) {};
\node[main_node] (12) at (\xD, \yC) {};
\node[main_node] (13) at (\xA, \yD) {$B_1$};
\node[main_node] (14) at (\xB, \yD) {$B_2$};
\node[main_node] (15) at (\xC, \yD) {$B_3$};
\node[main_node] (16) at (\xD, \yD) {$B_4$};

% 17 与 15 同高，18 与 0 同水平并与 4 同高
\node[small_node] (17) at (\xC, -1.3) {$B^*$};
\node[small_node] (18) at (1, 5.2) {};

% 绘制黑色连线
\path[draw, thick]
    (0) edge (1)
    (0) edge (2)
    (0) edge (3)
    (0) edge (4)
    (0) edge (18)
    (1) edge (2)
    (2) edge (6)
    (3) edge (7)
    (4) edge (8)
     (14) edge (15)
    (15) edge (16)
    (14) edge [bend left=23](16)
    (16) edge (17)
    (13) edge[ bend right=23] (15)
    (13) edge[ bend right=28] (16)
    (14) edge (17)
    (15) edge (17);

% 绘制红色连线
\path[draw, red, thick]
    (1) edge (5)
    (5) edge (9)
    (9) edge (13)
    (2) edge (6)
    (6) edge (10)
    (10) edge (14)
    (3) edge (7)
    (7) edge (11)
    (11) edge (15)
    (4) edge (8)
    (8) edge (12)
    (12) edge (16)
    (3) edge (18)
    (4) edge (18)
    (13) edge[ bend right=32] (17);
\end{tikzpicture}
\caption{A blow-up of $9$-framework with $m=1$ }
\label{fm_1}
\end{minipage}
\end{figure}
\begin{lemma}\label{2,3}
Let $G$ be an $\ell$-holed graph with $m=0$ and $2\leq k\leq 3$,
where for each $i\in[k]$, $s\lceil{\omega\over{\ell-1}}\rceil<|B_i|<(2s+2)\lceil{\omega\over{\ell-1}}\rceil$ and $|A_0|<(2s+2)\lceil{\omega\over{\ell-1}}\rceil$,
and $A_0$ and $B_i$ are colored with $\lceil{\ell\omega\over{\ell-1}}\rceil$ colors such that $c:=|(\cup^k_{i=1} c(B_i))\cap c(A_0)|\leq ks\lceil{\omega\over{\ell-1}}\rceil$.
Let $c_i:=|c(B_i)\cap c(A_0)|$ and 
when $k=2$, $c_2=\min \{c, s\lceil{\omega\over{\ell-1}}\rceil\}$, $c_1=c-c_2$; 
when $k=3$, $c_3=\min \{c, s\lceil{\omega\over{\ell-1}}\rceil\}$,
$c_2=\min \{c-c_3,s\lceil{\omega\over{\ell-1}}\rceil\}$,
$c_1=c-c_2-c_3$.
Let
$C_i:=c(B_i)\cap c(A_0)$ for each $i\in [k]$.
The vertices in $B_i$ with largest indices have colors $C_i$.
Let $V(G')=V(G)\cup V(A')$ where $A_0\cup A'$ is a clique and $G'[A_0\cup A', A_i]$ obeys the orderings of $A_0\cup A'$ and $A_i$ for each $i\in [k]$ and the ordering between $A_0$ and $A_i$ in $G'$ is the same as in $G$.
Then $G'$ is $\lceil{\ell\omega\over{\ell-1}}\rceil$-colorable and we can ensure that the common colors of $A_0$ and $B_i$ remain unchanged in $G'$
(See Figure \ref{f2,3}; A red line means that two sets obey the ordering, while a black line means two parts are complete.).
\end{lemma}

\pf By permuting the colors, we may assume $c(A_0)=\{\lceil{\ell\omega\over{\ell-1}}\rceil,\cdots, \lceil{\ell\omega\over{\ell-1}}\rceil-|A_0|+2\}$.
%Correspondingly, $C_i=c(B_i)\cap c(A_0)$ will also be mapped to $C'_i=c'(B_i)\cap c'(A_0)$ for each $i\in [k]$.
and $c(A')=\{\lceil{\ell\omega\over{\ell-1}}\rceil-|A_0|+1, \cdots, \lceil{\ell\omega\over{\ell-1}}\rceil-\omega+2\}$, so $c(A_0\cup A')=\{\lceil{\ell\omega\over{\ell-1}}\rceil,\cdots, \lceil{\ell\omega\over{\ell-1}}\rceil-\omega+2\}$.
Let $B'_i\subseteq B_i$ denote the order vertex set that has different color from $A_0$. %and we map $c(B'_i)$ to $c'(B'_i)$.

When $k=2$, $c(B_1)=\{2j+1:0\leq j\leq |B'_2|-1\}\cup \{2|B'_2|+1,\cdots, |B'_1|+|B'_2|\}\cup C_1$
and $c(B_2)=\{2j:1\leq j\leq |B'_2|\}\cup C_2$.
Let ${\cal{P}}_i$ have a balanced coloring.
When $|B'_2|\leq s\lceil{\omega\over{\ell-1}}\rceil$,
we define $L^c_{i,2s}=\{1,\cdots,\omega-1\}$.
When $|B'_2|>s\lceil{\omega\over{\ell-1}}\rceil$,
if $|B_2|\leq 2s\lceil{\omega\over{\ell-1}}\rceil$, we define $L_{2,2s,1}=\{1,\cdots, |B_2|\}$,
otherwise $L_{2,2s,1}=\{1,\cdots,2s\lceil{\omega\over{\ell-1}}\rceil\}\cup \{2s\lceil{\omega\over{\ell-1}}\rceil+2+2j:0\leq j\leq |B_2|-2s\lceil{\omega\over{\ell-1}}\rceil-1\}$;
if $|B_1|>s\lceil{\omega\over{\ell-1}}\rceil+|B'_2|$, we define 
$L_{1,2s,1}=\{1,\cdots,2s\lceil{\omega\over{\ell-1}}\rceil\}\cup \{2s\lceil{\omega\over{\ell-1}}\rceil+1+2j:0\leq j\leq |B'_2|-s\lceil{\omega\over{\ell-1}}\rceil-1\}\cup \{2|B'_2|+1,\cdots, |B_1|+|B'_2|-s\lceil{\omega\over{\ell-1}}\rceil\}$,
otherwise, $L_{1,2s,1}=\{1,\cdots,2s\lceil{\omega\over{\ell-1}}\rceil\}\cup \{2s\lceil{\omega\over{\ell-1}}\rceil+1+2j:0\leq j\leq |B_1|-2s\lceil{\omega\over{\ell-1}}\rceil-1\}$.
In all cases, $G'$ is $\lceil{\ell\omega\over{\ell-1}}\rceil$-colorable.

\medskip

When $k=3$, we have either $|B_2|<\lceil{3s\over 2}\lceil{\omega\over{\ell-1}}\rceil\rceil$ or $|B_2|\geq \lceil{3s\over 2}\lceil{\omega\over{\ell-1}}\rceil\rceil$.
We divide into cases. 

\medskip

\textbf{Case 1: $|B_2|<\lceil{3s\over 2}\lceil{\omega\over{\ell-1}}\rceil\rceil$}.\\
Let $i'\in \{0,1,2,3\}$ such that $|B'_{i'}|\geq \lceil{s\over 2}\lceil{\omega\over{\ell-1}}\rceil\rceil$ and $|B'_{i'+1}|< \lceil{s\over 2}\lceil{\omega\over{\ell-1}}\rceil\rceil$.
If for each $i\in [3]$, $|B'_{i}|\geq \lceil{s\over 2}\lceil{\omega\over{\ell-1}}\rceil\rceil$, then $i'=3$;
and if for each $i\in [3]$, $|B'_{i}|<\lceil{s\over 2}\lceil{\omega\over{\ell-1}}\rceil\rceil$, then $i'=0$.
When $i'\leq 1$, let $B'_1, B'_2, B'_3$ have a cyclic coloring and  $G'\backslash(A_0\cup A')$ have a balanced coloring, then 
$c(B_i)=c(B'_i)\cup C_i$ and
$L_{i, 2s, 1}=\{1,\cdots,|B_i|\}$ for each $i\in [3]$.
When $i'=2$, then $c_2<s\lceil{\omega\over{\ell-1}}\rceil$.
Let $B'_2=B'_{2,1}\cup B'_{2,2}$ where $|B'_{2,1}|=\lceil{s\over 2}\lceil{\omega\over{\ell-1}}\rceil\rceil$.
Then  let $B'_1, B'_{2,1}, B'_3$ and $B'_{2,2}$ have a cyclic coloring successively and  $G'\backslash(A_0\cup A')$ have a balanced coloring,
then $c(B_i)=c(B'_i)\cup C_i$
and
$L_{i, 2s, 1}=\{1,\cdots,|B_i|\}$ for each $i\in [3]$.
When $i'=3$, then $c_3<s\lceil{\omega\over{\ell-1}}\rceil$.
Let $B'_i=B'_{i,1}\cup B'_{i,2}$ where $|B'_{i,1}|=\lceil{s\over 2}\lceil{\omega\over{\ell-1}}\rceil\rceil$ for each $i\in \{2,3\}$.
Then  let $B'_1, B'_{2,1}, B'_{3,1}$ and $B'_{2,2}, B'_{3,2}$ have a cyclic coloring successively and  $G'\backslash(A_0\cup A')$ have a balanced coloring,
then $c(B_i)=c(B'_i)\cup C_i$
and
$L_{i, 2s, 1}=\{1,\cdots,|B_i|\}$ for each $i\in [3]$.

\medskip

\textbf{Case 2: $|B_2|\geq\lceil{3s\over 2}\lceil{\omega\over{\ell-1}}\rceil\rceil$}.\\
When $|B'_3|<\lceil{s\over 2}\lceil{\omega\over{\ell-1}}\rceil\rceil$,
let $B'_1, B'_2, B'_3$ have a cyclic coloring and  $G'\backslash(A_0\cup A')$ have a balanced coloring.
Otherwise, if $|B_3|\geq\lceil{3s\over 2}\lceil{\omega\over{\ell-1}}\rceil\rceil$, 
then $c_3=c<s\lceil{\omega\over{\ell-1}}\rceil$.
Let $B'_3=B'_{3,1}\cup B'_{3,2}$ where $|B'_{3,1}|=\lceil{s\over 2}\lceil{\omega\over{\ell-1}}\rceil\rceil$.
Let $B'_1, B'_2, B'_{3,1}$ and $ B'_{3,2}$ have a cyclic coloring successively and  $G'\backslash(A_0\cup A')$ have a balanced coloring.
Anyway, $G'$ is $\lceil{\ell\omega\over{\ell-1}}\rceil$-colorable.
\qed

\medskip

\begin{comment}

{\color{red} Difference between Cor 5.5 and Lem 3.8, Move after 3.8

%By Lemma \ref{2,3}, we can obtain the following corollary.

\begin{corollary}\label{cor}
Let $m\geq 1$ be an odd integer and $P := Z_0Z_1\cdots Z_m$ be the blow-up of a path where $Z_i$ is an ordered clique for each $i\in \{0,\cdots, m\}$ and
$G[Z_j, Z_{j+1}]$  obeys the orderings for each $j\in \{0,\cdots, m-1\}$.
Let $\omega := \omega(P)$ and $\chi'>\omega$ be an integer.
Suppose the vertices of $Z_0$ have colors $C_0$ in any order and the vertices of $Z_m$ have colors $C_m$ in any order such that $|C_0 \cap C_m|\leq {(m-1)(\chi'-\omega)\over 2}$ and $|C_0 \cup C_m|\leq \chi'$, then $P$ is $\chi'$-colorable.  
%{\color{red}It is worth noting that we do not have specific requirements for the coloring order of one of $Z_0$ and $Z_m$.}
\end{corollary}
}
\end{comment}

For the sake of convenience, we denote by $A_0$ the intersection of all directed paths in $T$ starting from $a_0$ to $a_i$ for each $i\in [k]$.

\begin{lemma}\label{m-neq-4}
$m\neq 4$.
\end{lemma}
\pf For otherwise, suppose $m=4$. We have $\ell=7$ and $s=1$, and since $k\geq 2m$, $k=9$ or $k=8$.
Let ${\cal P}^*_i=C_i$ for each $i\in [k]$.
Then $A_1$ is anticomplete to $A_{k-1}$.
Otherwise, by Lemma \ref{Claim 3.1.1}, $|A_0|\geq \lceil{\omega\over{\ell-1}}\rceil+1$, 
then $\sum^5_{i=0}|A_i|>\omega$, a contradiction.

Let $j$ be the largest index such that $A_1$ is complete to $A_j$ and anticomplete to $A_{j+1}$,
$5 \leq j \leq k-2$.
Let $A'$ be a vertex set such that $A'\cup A_0$ is a clique and $G[A_0\cup A', A_i]$ obeys the ordering for each $i\in \{j+1,\cdots,k\}$.
Then $G\backslash (\cup^k_{i=j+1} (A_{i}\cup C_{i})\cup A')$ is $\lceil{\ell\omega\over{\ell-1}}\rceil$-colorable.
Obviously, $|c(A_0)\cap c(\cup^k_{i=j+1}B_{i})|\leq (k-j)\lceil{\omega\over{\ell-1}}\rceil$.

When $k=9$, since we can color $A'$ such that
$|c(A_0\cup A')\cap c(\cup^9_{i=j+1}B_i)|\leq \sum^9_{i=j+1}|B_{i}|+|A_0|+|A'|-\lceil{\ell\omega\over{\ell-1}}\rceil<\omega-(j-4)\lceil{\omega\over{\ell-1}}\rceil+\omega-\lceil{\ell\omega\over{\ell-1}}\rceil<(9-j)\lceil{\omega\over{\ell-1}}\rceil$,
and adjust the colors in $\cup^9_{i=j+1}B_{i}$ such that for each $i\in \{j+1,\cdots,k\}$, $|c(B_{i}\cap c(A_0\cup A'))|\leq \lceil{\omega\over{\ell-1}}\rceil$.
By Lemma \ref{lem11}, $G$ is $\lceil{\ell\omega\over{\ell-1}}\rceil$-colorable.

When $k=8$, by Lemma \ref{2,3}, $G$ is $\lceil{\ell\omega\over{\ell-1}}\rceil$-colorable.
This completes the proof of the lemma.
\qed

\begin{lemma}
$m\neq 3$.
\end{lemma}
\pf For otherwise, suppose $m = 3$. So $6\leq k\leq 8$.

\medskip
\textbf{Case 1: $k=6$}\\
Then $A_1$ is complete to $A_5$.
Otherwise, let $A'$ be a vertex set such that $A'\cup A_0$ is a clique and $G[A_0\cup A', A_i]$ obeys the ordering for each $i\in \{5,6\}$.
Then $G\backslash (\cup^6_{i=5} (A_{i}\cup {\cal P}^*_{i})\cup A')$ is $\lceil{\ell\omega\over{\ell-1}}\rceil$-colorable.
Obviously, $|c(A_0)\cap c(\cup^6_{i=5}B_{i})|\leq 2\lceil{\omega\over{\ell-1}}\rceil$.
Adjust the colors in $\cup^6_{i=5}B_{i}$ such that for each $i\in \{5,6\}$, $|c(B_{i}\cap c(A_0\cup A'))|\leq \lceil{\omega\over{\ell-1}}\rceil$.
By Lemma \ref{2,3}, $G$ is $\lceil{\ell\omega\over{\ell-1}}\rceil$-colorable.
Let $B^*_1\subseteq V(S)\backslash V(\cup^6_{i=1}B_i)$ that is complete to $B_4, B_5, B_6$ and obeys the orderings with $B_1, B_2, B_3$. 
Since $A_1$ is complete to $A_5$, $|A_0|\geq s\lceil{\omega\over{\ell-1}}\rceil+1$.

If $A_2$ and $A_3$ are complete to $A_5$, 
let 
$B^*_2\subseteq V(S)\backslash V(\cup^6_{i=1}B_i\cup B^*_1)$ that is complete to $B_6$ and  $G[B_6\cup B^*_1\cup B^*_2, B_i]$ obeys the ordering for each $i\in [3]$.
Let $G\backslash(\cup^3_{i=1}({\cal P}^*_i\cup B_i)\cup A_6\cup {\cal P}^*_6\cup B^*_2)$ be colored 
by $\lceil{\ell\omega\over{\ell-1}}\rceil$ colors.
Since $|B_6\cup B^*_1|\leq 4s\lceil{\omega\over{\ell-1}}\rceil$,
then $|c(\cup^3_{i=0}A_i)\cap c(B_6\cup B^*_1)|\leq 4s\lceil{\omega\over{\ell-1}}\rceil$.
By adjusting the colors among $\cup^3_{i=0}A_i$, $|c(A_0)\cap c(B_6)|\leq|c(A_0)\cap c(B_6\cup B^*_1)|\leq s\lceil{\omega\over{\ell-1}}\rceil$ and $|c(A_1\cup A_2\cup A_3)\cap c(B_6\cup B^*_1)|\leq 3s\lceil{\omega\over{\ell-1}}\rceil$.
Then by Lemma \ref{lem11} and Lemma \ref{2,3}, $G$ is $\lceil{\ell\omega\over{\ell-1}}\rceil$-colorable.

If $A_2$ is complete to $A_5$, but $A_3$ is not complete to $A_5$,
let 
$B^*_2\subseteq V(S)\backslash V(\cup^6_{i=1}B_i\cup B^*_1)$ that is complete to $B_6$ and  $G[B_6\cup B^*_1\cup B^*_2, B_i]$ obeys the ordering for each $i\in [2]$ and 
$B^*_3\subseteq V(S)\backslash V(\cup^6_{i=1}B_i\cup B^*_1\cup B^*_2)$ that is complete to $B_4, B_5, B_6$ and  $G[B_6\cup B^*_1\cup B^*_3, B_3]$ obeys the ordering.
And let $A'_1\subseteq V(T)\backslash (\cup^6_{i=0}A_i)$
that is complete to $\cup^4_{i=0}A_i$ and $G[A_1\cup A_2\cup A'_1, A_5]$ obeys the ordering.
If $|B_6\cup B^*_1|\leq 3s\lceil{\omega\over{\ell-1}}\rceil$,
then $G\backslash({\cup^2_{i=1}({\cal P}^*_i\cup B_i)\cup A_6\cup {\cal P}^*_6}\cup B^*_2)$ is $\lceil{\ell\omega\over{\ell-1}}\rceil$-colorable.
By adjusting the colors among $\cup^2_{i=0}A_i$,
$|c(A_0)\cap c(B_6)|\leq|c(A_0)\cap c(B_6\cup B^*_1)|\leq s\lceil{\omega\over{\ell-1}}\rceil$ and $|c(A_1\cup A_2)\cap c(B_6\cup B^*_1)|\leq 2s\lceil{\omega\over{\ell-1}}\rceil$.
Then by Lemma \ref{lem11} and Lemma \ref{2,3}, $G$ is $\lceil{\ell\omega\over{\ell-1}}\rceil$-colorable.
Now, 
let $G\backslash ({\cup^3_{i=1}({\cal P}^*_i\cup B_i)\cup A_6\cup {\cal P}^*_6\cup A_5\cup {\cal P}^*_5}\cup B^*_2)$ be colored by $\lceil{\ell\omega\over{\ell-1}}\rceil$ colors.
Since $|B_5\cup B_6\cup B^*_1\cup B^*_3|<5s\lceil{\omega\over{\ell-1}}\rceil$,
by  adjusting the colors among $\cup^2_{i=0}A_i\cup A'_1$,
$|c(B_6)\cap c(A_0)|\leq s\lceil{\omega\over{\ell-1}}\rceil$,
$|c(B_6\cup B^*_1)\cap c(A_1\cup A_2)|\leq 2s\lceil{\omega\over{\ell-1}}\rceil$,
$|c(B_5\cup B_6\cup B^*_1\cup B^*_3)\cap c(A_3)|\leq s\lceil{\omega\over{\ell-1}}\rceil$,
and $|c(B_5)\cap c(A_0\cup A_1\cup A'_1)|\leq s\lceil{\omega\over{\ell-1}}\rceil$.
Then by Lemma \ref{lem11} and Lemma \ref{2,3}, $G$ is $\lceil{\ell\omega\over{\ell-1}}\rceil$-colorable.

If  $A_2$ and $A_3$ are not complete to $A_5$,
let $B^*_2, B^*_3\subseteq V(S)$ such that $B^*_2$ is complete to $B_4, B_5, B_6$ and obeys the orderings with $B_2, B_3$; and $B^*_3$ is complete to $ B_5, B_6$ and obeys the orderings with $B_2, B_3$.
Let $G\backslash ({\cup^3_{i=1}({\cal P}^*_i\cup B_i)\cup A_6\cup {\cal P}^*_6\cup A_5\cup {\cal P}^*_5}\cup B^*_3)$ be colored by $\lceil{\ell\omega\over{\ell-1}}\rceil$ colors.
Then
by  adjusting the colors among $\cup^2_{i=0}A_i\cup A'_1$,
$|c(B_6)\cap c(A_0)|\leq s\lceil{\omega\over{\ell-1}}\rceil$,
$|c(B_6\cup B^*_1)\cap c(A_1)|\leq s\lceil{\omega\over{\ell-1}}\rceil$,
$|c(B_5\cup B_6\cup B^*_1\cup B^*_2)\cap c(A_2\cup A_3)|\leq s\lceil{\omega\over{\ell-1}}\rceil$,
and $|c(B_5)\cap c(A_0\cup A_1\cup A'_1)|\leq s\lceil{\omega\over{\ell-1}}\rceil$.
Then by Lemma \ref{lem11} and Lemma \ref{2,3}, $G$ is $\lceil{\ell\omega\over{\ell-1}}\rceil$-colorable.

\medskip

\textbf{Case 2: $k=7$}\\
By Lemma \ref{2,3}, $A_1$ is complete to $A_6$.
Let $B^*_1\subseteq V(S)\backslash V(\cup^k_{i=1}B_i)$ that is complete to $\cup^k_{i=m+1}B_i$ and obeys the orderings with $B_1,\cdots, B_m$ respectively. 

If $A_2$ and  $A_3$ are complete to $A_6$,
let $B^*_2\subseteq V(S)$ that is complete to $B_7$ and obeys the orderings with $B_1,B_2, B_3$ respectively. 
Then $|c(B_7\cup B^*_1)\cap c(A_1\cup A_2\cup A_3)|\leq 3s\lceil{\ell\omega\over{\ell-1}}\rceil$.
By adjusting the colors among $\cup^3_{i=1}A_i$,
and by Lemma \ref{2,3}, $G$ is $\lceil{\ell\omega\over{\ell-1}}\rceil$-colorable.

If $A_2$ is complete to $A_6$, but $A_3$ is not complete to $A_6$,
let $B^*_2\subseteq V(S)$ that is complete to $B_7$ and obeys the orderings with $B_1,B_2$ respectively. 
Then $|c(B_7\cup B^*_1)\cap c(A_1\cup A_2\cup A_0)|\leq 3s\lceil{\ell\omega\over{\ell-1}}\rceil$.
By adjusting the colors among $\cup^2_{i=0}A_i$,
and by Lemma \ref{2,3}, $G$ is $\lceil{\ell\omega\over{\ell-1}}\rceil$-colorable.

If $A_2$ and  $A_3$ are complete to $A_5$, and are anticomplete to $A_6$,
let $B^*_2$ be complete to $\cup^7_{i=4}B_i$ and obey the orderings with $B_2, B_3$ respectively.
Let $B^*_3$ be complete to $\cup^7_{i=6}B_i$ and obey the orderings with $B_2, B_3$ respectively.
Let $A'_1\subseteq V(T)$ be complete to $\cup^5_{i=0}A_i$ and obey the ordering with $A_6$.
Let $G\backslash (\cup^7_{i=6}({\cal P}^*_i\cup A_i)\cup B_1\cup {\cal P}^*_1)$ be colored by  be colored by $\lceil{\ell\omega\over{\ell-1}}\rceil$ colors.
Then $|c(B_6\cup B_7\cup B^*_1)\cap c(A_0\cup A_1\cup A'_1)|<3s\lceil{\omega\over{\ell-1}}\rceil$.
By adjusting the colors among $\cup^1_{i=0}A_i\cup A'_1$,
$|c(B_7)\cap c(A_0)|\leq s\lceil{\omega\over{\ell-1}}\rceil$,
$|c(B_7\cup B^*_1)\cap c(A_1)|\leq s\lceil{\omega\over{\ell-1}}\rceil$,
$|c(B_6)\cap c(A_0\cup A_1\cup A'_1)|\leq s\lceil{\omega\over{\ell-1}}\rceil$.
Then by Lemma \ref{lem11}, $G$ is $\lceil{\ell\omega\over{\ell-1}}\rceil$-colorable.

If $A_2$ is complete to $A_5$ and $A_3$ is anticomplete to $A_5$, 
let $B^*_2\subseteq V(S)$ be complete to $\cup^7_{i=4}B_i$ and obey the orderings with $B_2, B_3$ respectively.
Let $A'_1\subseteq V(T)$ be complete to $\cup^5_{i=0}A_i$ and obey the ordering with $A_6$.
Then $|B^*_1\cup B_7|>2s\lceil{\ell\omega\over{\ell-1}}\rceil$ and $|B^*_1\cup B_7\cup B_6\cup B^*_2|<4s\lceil{\ell\omega\over{\ell-1}}\rceil$.
Let $G\backslash (\cup^2_{i=1}({\cal P}^*_i\cup B_i)\cup (\cup^7_{i=6}A_i\cup {\cal P}^*_i))$ be colored by  be colored by $\lceil{\ell\omega\over{\ell-1}}\rceil$ colors.
Then
by  adjusting the colors among $\cup^3_{i=0}A_i$,
$|c(B_7)\cap c(A_0)|\leq s\lceil{\omega\over{\ell-1}}\rceil$,
$|c(B_7\cup B^*_1)\cap c(A_1)|\leq s\lceil{\omega\over{\ell-1}}\rceil$,
$|c(B_6\cup B_7\cup B^*_1\cup B^*_2)\cap c(A_2)|\leq s\lceil{\omega\over{\ell-1}}\rceil$
and $|c(B_6)\cap c(A_0\cup A_1\cup A'_1)|\leq s\lceil{\omega\over{\ell-1}}\rceil$.
Then by Lemma \ref{lem11}, $G$ is $\lceil{\ell\omega\over{\ell-1}}\rceil$-colorable.

If $A_2$ and  $A_3$ are anticomplete to $A_5$, 
let $A'_1\subseteq V(T)$ such that $A'_1$ is complete to $\cup^4_{i=1}A_i$ and obeys the orderings with $A_5, A_6$ respectively.
Let $A'_2\subseteq V(T)$ such that $A'_1$ is complete to $A_0\cup A_1\cup A'_1$ and obeys the orderings with $A_5, A_6$ respectively.
Then $|B^*_1\cup B_7|>2s\lceil{\ell\omega\over{\ell-1}}\rceil$ and $|c(B^*_1\cup B_7\cup B_6\cup B_5)\cap c(A_0\cup A_1\cup A'_1)|<3s\lceil{\ell\omega\over{\ell-1}}\rceil$.
Let $G\backslash (\cup^7_{i=5}({\cal P}^*_i\cup A_i)\cup B_1\cup {\cal P}^*_1)$ be colored by $\lceil{\ell\omega\over{\ell-1}}\rceil$ colors.
Then by adjusting the colors among $\cup^1_{i=0}A_i\cup A'_1$,
$|c(B_7)\cap c(A_0)|\leq s\lceil{\omega\over{\ell-1}}\rceil$,
$|c(B_7\cup B^*_1)\cap c(A_1)|\leq s\lceil{\omega\over{\ell-1}}\rceil$
and 
$|c(B_6\cup B_5)\cap c(A_0\cup A_1\cup A'_1)|\leq 2s\lceil{\omega\over{\ell-1}}\rceil$.
Then by Lemma \ref{lem11} and Lemma \ref{2,3}, $G$ is $\lceil{\ell\omega\over{\ell-1}}\rceil$-colorable.

\medskip

\textbf{Case 3: $k=8$}\\
By Lemma \ref{2,3}, $A_1$ is complete to $A_7$ or $A_1$ is anticomplete to $\cup^8_{i=5}A_i$.

If $A_1$ is anticomplete to $\cup^8_{i=5}A_i$,
let $A'\subseteq V(T)$ such that $A'$ is complete to $A_0$ and $G[A_0\cup A', A_i]$ obey the orderings for each $i\in \{5,6,7,8\}$.
Let $G\backslash \cup^8_{i=5}(A_i\cup {\cal P}^*_i)\cup A'$ be colored by $\lceil{\ell\omega\over{\ell-1}}\rceil$ colors.
Then $|c(A_0)\cap c(\cup^8_{i=5}B_i)|\leq 4s\lceil{\omega\over{\ell-1}}\rceil$.
And $|c(A_0\cup A')\cap c(\cup^8_{i=5}B_i)|\leq 4s\lceil{\omega\over{\ell-1}}\rceil$,
 $G$ is $\lceil{\ell\omega\over{\ell-1}}\rceil$-colorable.

If $A_1$ is complete to $A_7$, let $B^*_1\subseteq V(S)\backslash V(\cup^8_{i=1}B_i)$ that is complete to $\cup^8_{i=4}B_i$ and obeys the ordering with $B_1$ respectively. 
Then $|B_8\cup B^*_1|>2s\lceil{\omega\over{\ell-1}}\rceil$ and $|\cup^8_{i=4}B_i\cup B^*_1|>\omega$, a contradiction.
\qed

\begin{lemma}
$m\neq 2$
\end{lemma}
\pf Suppose not, then we have $4\leq k\leq 7$.
The proof when $k=4$ is similar to the case when $m=0$ and $k=4$.
So we consider the case when $5\leq k\leq 7$.

\medskip 
\textbf{Case 1: $k=5$}\\
By Lemma \ref{2,3}, $A_1$ is complete to $A_4$.
If $A_2$ is complete to $A_4$,
then let $A'_1\subseteq V(T)$ such that $A'_1$ is complete to $\cup^2_{i=1}A_i$ and obeys the orderings with $A_3$, $A_4$.
Let $B^*_1, B^*_2\subseteq V(S)$ such that $B^*_1$ is complete to $\cup^5_{i=3}B_i$, $B^*_2$ is complete to $B_5$ and  $G[B_5\cup B^*_1\cup B^*_2, B_i]$ obeys the orderings for each $i\in \{2,3\}$.
Let $A_0\cup A_1\cup A_2\cup A'_1\cup B_3\cup B_4\cup B_5\cup B^*_1$ be colored by $\lceil{\ell\omega\over{\ell-1}}\rceil$ colors.
Then $|c(B_3\cup B_4\cup B_5\cup B^*_1)\cap c(A_0\cup A_1\cup A_2\cup A'_1)|\leq 5s\lceil{\ell\omega\over{\ell-1}}\rceil$.
By adjusting the colors among $A_0\cup A_1\cup A_2\cup A'_1$,
$|c(B_3\cup B_4)\cap c(A_0\cup A_1\cup A_2\cup A'_1)|\leq 2s\lceil{\ell\omega\over{\ell-1}}\rceil$,
we have $|c( B_5)\cap c(A_0)|\leq s\lceil{\ell\omega\over{\ell-1}}\rceil$
and $|c(B_5\cup B^*_1)\cap c( A_1\cup A_2)|\leq 2s\lceil{\ell\omega\over{\ell-1}}\rceil$.
Then by Lemmas \ref{lem11} and \ref{2,3}, $G$ is $\lceil{\ell\omega\over{\ell-1}}\rceil$-colorable.

If $A_2$ is anticomplete to $A_4$, then
let $A'_1\subseteq V(T)$ such that $A'_1$ is complete to $\cup^3_{i=1}A_i$ and obeys the ordering with $A_4$.
Let $B^*_1, B^*_2\subseteq V(S)$ such that $B^*_1$ is complete to $\cup^5_{i=2}B_i$ and  obeys the ordering with $B_1$
while $B^*_2$ is complete to $\cup^5_{i=3}B_i$ and  obeys the ordering with $ B_2$.

Let $\cup^3_{i=0}A_i\cup A'_1\cup B_4\cup B_5\cup B^*_1\cup B^*_2$ be colored by $\lceil{\ell\omega\over{\ell-1}}\rceil$ colors.
Then $|c(B_4\cup B_5\cup B^*_1\cup B^*_2)\cap c(\cup^3_{i=0}A_i\cup A'_1)|\leq 5s\lceil{\ell\omega\over{\ell-1}}\rceil$.
By adjusting the colors, we have
$|c(B_4)\cap c(A_0\cup A_1\cup  A'_1)|\leq s\lceil{\ell\omega\over{\ell-1}}\rceil$,
$|c(B_4\cup B_5\cup B^*_1\cup B^*_2)\cap c(A_i)|\leq s\lceil{\ell\omega\over{\ell-1}}\rceil$ for each $i\in \{2,3\}$
and $|c(B_5\cup B^*_1)\cap c( A_1\cup A_2)|\leq 2s\lceil{\ell\omega\over{\ell-1}}\rceil$.
Then by Lemmas \ref{lem11} and \ref{2,3}, $G$ is $\lceil{\ell\omega\over{\ell-1}}\rceil$-colorable.

\medskip 

\textbf{Case 2: $k=6$}\\
By Lemma \ref{2,3}, $A_1$ is complete to $A_5$,
then $|A_0|\geq s\lceil{\omega\over{\ell-1}}\rceil+1$.

If $A_2$ is complete to $A_5$, then
let  $B^*_1,B^*_2\subseteq V(S)\backslash V(\cup^6_{i=1}B_i)$ such that $B^*_1$ is complete to $\cup^6_{i=3}B_i$ and obeys the orderings with $B_1, B_2$ and $B^*_2$ is complete to $B_6$ and obeys the orderings with $B_1, B_2$.
Let $G\backslash \cup^2_{i=1}(B_i\cup {\cal P}^*_i)\cup A_6\cup {\cal P}^*_6\cup B^*_2$ be colored by  $\lceil{\ell\omega\over{\ell-1}}\rceil$ colors.
Then $|c(\cup^2_{i=0}A_i)\cap c(B_6\cup B^*_1)|\leq 3s\lceil{\omega\over{\ell-1}}\rceil$.
By adjusting the colors among $\cup^2_{i=0}A_i$,
we have $|c(A_0)\cap c(B_6)|\leq s\lceil{\omega\over{\ell-1}}\rceil$
and $|c(A_1\cup A_2)\cap c(B_6\cup B^*_1)|\leq 2s\lceil{\omega\over{\ell-1}}\rceil$.
By Lemmas \ref{2,3} and \ref{lem11},
 $G$ is $\lceil{\ell\omega\over{\ell-1}}\rceil$-colorable.

 If $A_2$ is anticomplete to $A_5$, then
 let $B^*_1\subseteq V(S)\backslash V(\cup^6_{i=1}B_i)$ such that $B^*_1$ is complete to $\cup^6_{i=3}B_i$ and obeys the ordering with $B_1$.
Then $|B_6\cup B^*_1|\geq 2s\lceil{\omega\over{\ell-1}}\rceil$.
If $A_2 $ is anticomplete to $A_4$, then
let $A'_1, A'_2\subseteq V(T)\backslash V(\cup^6_{i=1}A_i)$ such that $A'_1$ is complete to $\cup^3_{i=1}A_i$, $A'_2$ is complete to $A_1$ and $G[A_1\cup A'_1\cup A'_2, A_i]$ obey the orderings for each $i\in \{4,5\}$.
Let $G\backslash \cup^6_{i=4}(A_i\cup {\cal P}^*_i)\cup B_1\cup {\cal P}^*_1\cup A'_2$ be colored by  $\lceil{\ell\omega\over{\ell-1}}\rceil$ colors.
Then $|c(\cup^1_{i=0}A_i\cup A'_1)\cap c(\cup^6_{i=4}B_i\cup B^*_1)|\leq 4s\lceil{\omega\over{\ell-1}}\rceil$.
By adjusting the colors among $\cup^6_{i=4}B_i\cup B^*_1$,
we have $|c(A_0)\cap c(B_6)|\leq s\lceil{\omega\over{\ell-1}}\rceil$, 
$|c(A_1)\cap c(B_6\cup B^*_1)|\leq s\lceil{\omega\over{\ell-1}}\rceil$
and $|c(A_1\cup A_0\cup A'_1)\cap c(B_4\cup B_5)|\leq 2s\lceil{\omega\over{\ell-1}}\rceil$.
By Lemmas \ref{2,3} and \ref{lem11},
 $G$ is $\lceil{\ell\omega\over{\ell-1}}\rceil$-colorable.
If $A_2 $ is complete to $A_4$, then
 let $B^*_2\subset V(S)\backslash V(\cup^6_{i=1}B_i)$ such that $B^*_2$ is complete to $\cup^6_{i=3}B_i$ and $G[B_5\cup B_6\cup B^*_1\cup B^*_2, B_2]$ obeys the ordering.
 Let $A'_1\subseteq V(T)\backslash V(\cup^6_{i=1}A_i)$ such that $A'_2$ is complete to $\cup^4_{i=1}A_i$ and obeys the ordering with $A_5$.
 Let $G\backslash \cup^6_{i=5}(A_i\cup {\cal P}^*_i)\cup B_1\cup B_2\cup (\cup^2_{i=1} {\cal P}^*_i)$ be colored by  $\lceil{\ell\omega\over{\ell-1}}\rceil$ colors.
Then $|c(\cup^2_{i=0}A_i\cup A'_1)\cap c(\cup^6_{i=5}B_i\cup B^*_1\cup B^*_2)|\leq 4s\lceil{\omega\over{\ell-1}}\rceil$.
By adjusting the colors $\cup^2_{i=0}A_i\cup A'_1$,
we have $|c(A_0)\cap c(B_6)|\leq s\lceil{\omega\over{\ell-1}}\rceil$, 
$|c(A_1)\cap c(B_6\cup B^*_1)|\leq s\lceil{\omega\over{\ell-1}}\rceil$,
$|c(A_2)\cap c(B_5\cup B_6\cup B^*_1\cup B^*_2)|\leq s\lceil{\omega\over{\ell-1}}\rceil$
and $|c(A_1\cup A_0\cup A'_1)\cap c( B_5)|\leq s\lceil{\omega\over{\ell-1}}\rceil$.
By Lemma \ref{lem11},
 $G$ is $\lceil{\ell\omega\over{\ell-1}}\rceil$-colorable.

\medskip

\textbf{Case 3: $k=7$}\\
By Lemma \ref{2,3}, $A_1$ is complete to $A_6$ or $A_1$ is anticomplete to $\cup^7_{i=4}A_i$.

If $A_1$ is anticomplete to $\cup^7_{i=4}A_i$, then
let $A'\subseteq V(T)$ such that $A'$ is complete to $A_0$ and $G[A_0\cup A', A_i]$ obey the orderings for each $i\in \{4,5,6,7\}$.
Let $G\backslash \cup^7_{i=4}(A_i\cup {\cal P}^*_i)\cup A'$ be colored by  $\lceil{\ell\omega\over{\ell-1}}\rceil$ colors.
Then $|c(A_0)\cap c(\cup^7_{i=4}B_i)|\leq 4s\lceil{\omega\over{\ell-1}}\rceil$ and $|c(A_0\cup A')\cap c(\cup^7_{i=4}B_i)|\leq 4s\lceil{\omega\over{\ell-1}}\rceil$.
Rhen for each $i\in \{4,\cdots,7\}$, we have $|c(A_0\cup A')\cap c(B_i)|\leq s\lceil{\omega\over{\ell-1}}\rceil$.
By  Lemma \ref{lem11},
 $G$ is $\lceil{\ell\omega\over{\ell-1}}\rceil$-colorable.

If $A_1$ is complete to $A_6$, let $B^*_1\subseteq V(S)\backslash V(\cup^7_{i=1}B_i)$ that is complete to $\cup^7_{i=3}B_i$ and obeys the ordering with $B_1$. 
Either by Lemma \ref{2,3},  $G$ is $\lceil{\ell\omega\over{\ell-1}}\rceil$-colorable;
or $|B_7\cup B^*_1|>2s\lceil{\omega\over{\ell-1}}\rceil$ and $|\cup^7_{i=3}B_i\cup B^*_1|>\omega$, a contradiction.
\qed

\begin{lemma}
$m\neq 1$.
\end{lemma}
\pf
Suppose not, then $3 \leq k\leq 6$.
%It is impossible that $k=2$.
The proof of the case when $k=3$ is similar to that of the case when $m=0$ and $k=3$.
So we assume that $4\leq k\leq 6$.

\medskip

\textbf{Case 1: $k=6$}\\
By Lemma \ref{2,3}, $A_1$ is complete to $A_5$ or $A_1$ is anticomplete to $\cup^6_{i=3}A_i$.

If $A_1$ is complete to $A_5$, then let $B^*_1\subseteq  V(S)\backslash V(\cup^6_{i=1}B_i)$ such that $B^*_1$ is complete to $\cup^6_{i=3}B_i$ and obeys the ordering with $B_1$. 
If $|B_6\cup B^*_1|>2s\lceil{\omega\over{\ell-1}}\rceil$, then $|\cup^6_{i=2}B_i\cup B^*_1|>\omega$, a contradiction.
Otherwise, by Lemma \ref{2,3},  $G$ is $\lceil{\ell\omega\over{\ell-1}}\rceil$-colorable.

If $A_1$ is anticomplete to $\cup^6_{i=3}A_i$, then let $A'\subseteq V(T)$ such that $A'$ is complete to $A_0$ and $G[A_0\cup A', A_i]$ obey the orderings for each $i\in \{3,4,5,6\}$.
Let $G\backslash \cup^6_{i=3}(A_i\cup {\cal P}^*_i)\cup A'$ be colored by  $\lceil{\ell\omega\over{\ell-1}}\rceil$ colors.
Then we have $|c(A_0)\cap c(\cup^6_{i=3}B_i)|\leq 4s\lceil{\omega\over{\ell-1}}\rceil$ and $|c(A_0\cup A')\cap c(\cup^6_{i=3}B_i)|\leq 4s\lceil{\omega\over{\ell-1}}\rceil$.
Then for each $i\in \{3,\cdots,6\}$, $|c(A_0\cup A')\cap c(B_i)|\leq s\lceil{\omega\over{\ell-1}}\rceil$.
By Lemma \ref{lem11},
 $G$ is $\lceil{\ell\omega\over{\ell-1}}\rceil$-colorable.

\medskip

\textbf{Case 2: $k=5$}\\
By Lemma \ref{2,3},
$A_1$ is complete to $A_4$.
Let $B^*_1\subseteq V(S)\backslash V(\cup^5_{i=1}B_i)$ such that $B^*_1$ is complete to $\cup^5_{i=2}B_i$ and obeys the ordering with $B_1$.
Let $A'_1\subseteq V(T)\backslash V(\cup^5_{i=1}A_i)$ such that $G[A_1\cup A'_1, A_i]$  obeys the orderings for each $i\in \{2,3,4\}$.
Then $|A_0|\geq s\lceil{\omega\over{\ell-1}}\rceil+1$
and $|B_5\cup B^*_1|\leq 2s\lceil{\omega\over{\ell-1}}\rceil$.
Let $\cup^5_{i=2}B_i\cup A_0\cup A_1\cup A'_1\cup B^*_1$ be colored by $\lceil{\ell\omega\over{\ell-1}}\rceil$ colors.
Then $|c(A_0\cup A_1\cup A'_1)\cap c(\cup^5_{i=2}B_i\cup B^*_1)|\leq 5s\lceil{\omega\over{\ell-1}}\rceil$.
By adjusting the colors among $A_0\cup A_1\cup A'_1$,
we have $|c(A_0)\cap c(B_5)|\leq s\lceil{\omega\over{\ell-1}}\rceil$,
$|c(A_1)\cap c(B_5\cup B^*_1)|\leq s\lceil{\omega\over{\ell-1}}\rceil$,
and $|c(A_0\cup A_1\cup A'_1)\cap c(B_i)|\leq s\lceil{\omega\over{\ell-1}}\rceil$ for each $i\in\{2,3,4\}$.
By Lemma \ref{lem11}, $G$ is $\lceil{\ell\omega\over{\ell-1}}\rceil$-colorable.

\medskip

\textbf{Case 3: $k=4$}\\
By Lemma \ref{2,3},
$A_1$ is complete to $A_3$.
Let $B^*_1\subseteq V(S)\backslash V(\cup^4_{i=1}B_i)$ such that $B^*_1$ is complete to $\cup^4_{i=2}B_i$ and obeys the ordering with $B_1$.
Then $|A_0| \geq s\lceil{\omega\over{\ell-1}}\rceil+1$ and $|B_4\cup B^*_1| \leq 2s\lceil{\omega\over{\ell-1}}\rceil$.
Let $\cup^4_{i=2}B_i\cup A_0\cup A_1\cup B^*_1$ be colored by $\lceil{\ell\omega\over{\ell-1}}\rceil$ colors.
Now suppose $|c(A_0\cup A_1)\cap c(\cup^4_{i=2}B_i\cup B^*_1)|\leq 4s\lceil{\omega\over{\ell-1}}\rceil$.
By adjusting the colors among $\cup^4_{i=2}B_i\cup B^*_1$,
we have $|c(A_0)\cap c(B_4)|\leq s\lceil{\omega\over{\ell-1}}\rceil$,
$|c(A_1)\cap c(B_4\cup B^*_1)|\leq s\lceil{\omega\over{\ell-1}}\rceil$
and $|c(A_0\cup A_1)\cap c(B_2\cup B_3)|\leq 2s\lceil{\omega\over{\ell-1}}\rceil$.
By Lemmas \ref{lem11} and \ref{2,3}, $G$ is $\lceil{\ell\omega\over{\ell-1}}\rceil$-colorable.
So
$4s\lceil{\omega\over{\ell-1}}\rceil<|c(A_0\cup A_1)\cap c(\cup^4_{i=2}B_i\cup B^*_1)|<5s\lceil{\omega\over{\ell-1}}\rceil$.
Let $C$ be the set of $\lceil{\ell\omega\over{\ell-1}}\rceil$ colors, $Y=C\backslash c(A_0\cup A_1)$ 
and $X=c(A_0\cup A_1)\cap c(B_2\cup B_3\cup B_4\cup B^*_1)$.
Then $|Y|>(s+1)\lceil {\omega \over {\ell-1}}\rceil$
and $|X|>4s\lceil{\omega\over{\ell-1}}\rceil$.
Let $X'$, $Y'$ be the subset of $X$ and $Y$, respectively, where $|X'|=|X|-4s\lceil {\omega \over {\ell-1}}\rceil$.
If there exists $Y'\cap {c(B_2\cup B_3\cup B_4\cup B^*_1)}=\emptyset$ such that $|c(A_0\cup A_1)\backslash X'\cup Y'|=|c(A_0\cup A_1)|$,
then by adjusting the colors, $G$ is $\lceil{\ell\omega\over{\ell-1}}\rceil$-colorable, a contradiction.
So for any $Y'$, $|Y'|<|X'|$, that is, $Y\backslash Y'\subseteq (c(B_2\cup B_3\cup B_4\cup B^*_1)\backslash X)$.
Then $|c(B_2\cup B_3\cup B_4\cup B^*_1)|>(s+1)\lceil{\omega\over{\ell-1}}\rceil-|X|+4s\lceil{\ell\omega\over{\ell-1}}\rceil+|X|>\omega$, a contradiction.\qed

\medskip

Therefore, by the above lemmas, when $\ell \equiv 3 \pmod 4$, the blow-up of $\ell$-frameworks $G$ is $\lceil{\ell\over{\ell-1}}\omega(G)\rceil$-colorable.

\subsection{The case when $\ell \equiv 1 \pmod 4$}

Let $\ell=4s+1$ with $s\geq 2$
and $G$ be a minimal counterexample of Lemma \ref{lem-color-framework} with minimum $|V(G)|$.
Recall that both $\bigcup^{m+1}_{i=1}A_i$ and $\bigcup^{k}_{j=m+1}B_j$ are cliques.
Then by Lemma \ref{Claim 3.1.1}, $m$ is at most $2$ and $k-m$ is at most $3$.
%Let $A^*$ denote the set of common vertices of the path from $a_0$ to $a_i$ for every $i\in [k]$ in $T$.

\begin{lemma}
    $m\neq 0$.
\end{lemma}
\pf Suppose for the contradiction that $m=0$.
Since the blow-up of an $\ell$-cycle  is ${\lceil {\ell\omega \over {\ell-1}}\rceil}$-colorable by Lemma \ref{lem:blow-up-cycle-colorable}, we may assume that $k=3$.

Let $m_i:=|B_i|$ for each $i\in [3]$.
We may assume that $m_1\geq m_2\geq m_3\geq {\lceil {\omega \over {\ell-1}}\rceil}s+1$.
And we assume that $|A_0|=\omega-1$, let $A_0$ and $A_i$ obey the orderings for each $i\in [3]$.
Let the coloring of $A_0$ be $A^c_0=\{1,\cdots, \omega-1 \}$.
Let $\alpha=\lceil{3(s-1)\over 2} {\lceil {\omega \over {\ell-1}}\rceil} \rceil $ and $\beta=\lceil{3s\over 2} {\lceil {\omega \over {\ell-1}}\rceil} \rceil $.

If $m_1\leq 1+{3s \over 2} \lceil{\omega\over{\ell-1}}\rceil$,
then $(B_1, B_2, B_3)$ have a cyclic coloring and $G\backslash A_0$ have a balanced coloring with regard to such cyclic coloring of $(B_1, B_2, B_3)$.
Let $B_{i,1}=\{3j+i: 0 \le j \le m_3-1\}$ for $i \in [3]$
%Let $B_{1,1}=\{1, \cdots, 3m_3-2\}$, $B_{2,1}=\{2, \cdots, 3m_3-1\}$, $B_{3,1}=\{3, \cdots, 3m_3\}$ and the elements in $B_{i,1}$ for each $i\in [3]$ form an arithmetic sequence with the first term $i$ and common difference $3$.
 $B_{i,2}=\{3m_3+i+2j: 0 \le j \le m_2 - m_3-1\}$ for $i \in [2]$, and $B_{1,3}=\{2m_2+m_3+1, \cdots, m_1+m_2+m_3\}$.
%Let $B_{1,2}=\{3m_3+1, \cdots, 2m_2+m_3-1\}$, $B_{2,2}=\{3m_3+2, \cdots, 2m_2+m_3\}$ and the elements in $B_{i,2}$ for each $i\in [2]$ form an arithmetic sequence with the first term $3m_3+i$ and common difference $2$.
%Let $B_{1,3}=\{2m_2+m_3+1, \cdots, m_1+m_2+m_3\}$  $B_{1,4}=\{3m_3+1, \cdots, m_1+m_2+m_3\}$.
%and the elements in $B_{1,3}$ and $B_{1,4}$ form an arithmetic sequence with the first term $2m_2+m_3+1$, $3m_3+1$ and common difference $1$, respectively.
Note that all sets are ordered, and some may be empty.
Hence,  $(B_1, B_2, B_3)$ has a cyclic coloring such that the vertices in $B_1$, $B_2$ and $B_3$ are colored with $B_{1,1}\cup B_{1,2}\cup B_{1,3}$, $B_{2,1}\cup B_{2,2}$ and $B_{3,1}$, respectively.

If $m_1> 1+{3s \over 2} \lceil{\omega\over{\ell-1}}\rceil$,
then $(B_1, B'_2, B'_3)$ have a cyclic coloring and
$(B_2\backslash B'_2$, $B_3\backslash B'_3)$ have a cyclic coloring where $|B'_2|={{s\over 2}\lceil{\omega\over{\ell-1}}\rceil+{g\over2}}$ and $|B'_3|={{s\over 2}\lceil{\omega\over{\ell-1}}\rceil-{g\over2}}$ with $g=1$ if $\beta\equiv 2\pmod{3}$, and $g=0$ otherwise.
For $i\in [2]$, 
let $B_{i,1}=\{3j+i:0\leq j\leq |B'_2|-1\}$
and $B_{3,1}=\{3j:1\leq j\leq |B'_3|\}$.
The order of the elements is from left bound of $j$ to the right bound of $j$.
%Let $B_{1,1}=\{1, \cdots,{\lceil{{3s\lceil{\omega\over{\ell-1}}\rceil}\over 2}\rceil-2}+{g}\}$, $B_{2,1}=\{2, \cdots, \lceil{{3s\lceil{\omega\over{\ell-1}}\rceil}\over 2}\rceil-1+{g}\}$, $B_{3,1}=\{3, \cdots, \lceil{{3s\lceil{\omega\over{\ell-1}}\rceil}\over 2}\rceil-{2g}\}$ 
%and the elements in $B_{i,1}$ for each $i\in [3]$ form an arithmetic sequence with the first term $i$ and common difference $3$.
Let $B_{1,2}=\{\beta+1, \cdots, m_1+s\lceil{\omega\over{\ell-1}}\rceil\}$,
$B_{i,2}=\{2m_3+m_1+i-3-2j: m_3-|B'_i|-1\geq j\geq 0 \}$ for $i\in \{2,3\}$,
%$B_{i,2}=\{m_1+s\lceil{\omega\over{\ell-1}}\rceil+i-1+(-1)^ig,\cdots,2m_3+m_1+i-3 \}$
%$B_{2,2}=\{m_1+s\lceil{\omega\over{\ell-1}}\rceil+1+{g}, \cdots, 2m_3+m_1-1\}$ and
%$B_{3,2}=\{m_1+s\lceil{\omega\over{\ell-1}}\rceil+2-{g}, \cdots, 2m_3+m_1\}$
%and the elements in $B_{1,2}$ form an arithmetic sequence with the first term $\lceil{{3s\lceil{\omega\over{\ell-1}}\rceil}\over 2}\rceil+1$ and common difference $1$, the elements in $B_{i,2}$ for each $i\in \{2, 3\}$ form an arithmetic sequence with the common difference $2$.
and $B_{2,3}=\{2m_3+m_1+1, \cdots, m_1+m_2+m_3\}$.
%and the elements in $B_{2,3}$ form an arithmetic sequence with the first term $2m_3+m_1+1$ and common difference $1$.
Note that all sets are ordered, and some may be empty.
Hence,   $B_1$, $B_2$ and $B_3$ are colored with $B_{1,1}\cup B_{1,2}$, $B_{2,1}\cup B_{2,2}\cup B_{2,3}$ and $B_{3,1}\cup B_{3,2}$, respectively.

%\begin{enumerate}[i)]
%    \item  If $m_2\neq m_3$, then the vertices in $B_1$, $B_2$ and $B_3$ are colored one-to-one with the elements in $B_{1,1}\cup B_{1,2}$, $B_{2,1}\cup B_{2,2}\cup B_{2,3}$ and $B_{3,1}\cup B_{3, 2}$, respectively, following the order of their indices from smallest to largest.
%    \item  If $m_2=m_3$, then the vertices in in $B_1$, $B_2$ and $B_3$ are colored one-to-one with the elements in $B_{1,1}\cup B_{1,2}$, $B_{2,1}\cup B_{2,2}$ and $B_{3,1}\cup B_{3, 2}$, respectively, following the order of their indices from smallest to largest.
%\end{enumerate}

We define
$X_1=\{1, \cdots, \alpha-2, \alpha-1, \alpha\},$
$X_2=\{1, \cdots, \alpha-2, \alpha, \alpha+1\},$ and
$X_3=\{1, \cdots, \alpha-2, \alpha-1, \alpha+1\}.$
Now we verify that the balanced coloring gives a proper coloring of $G$.
We divide into several cases according to the values of $m_1$ and residue of $\alpha$ modulo $3$.

\begin{comment}
\begin{enumerate}[i)]
    \item  \label{2.1}If $m_1\neq m_2$ and $m_2\neq m_3$, then the vertices in $B_1$, $B_2$ and $B_3$ are colored one-to-one with the elements in $B_{1,1}\cup B_{1,2}\cup B_{1,3}$, $B_{2,1}\cup B_{2,2}$ and $B_{3,1}$, respectively, following the order of their indices from smallest to largest.
    \item \label{2.2} If $m_1=m_2=m_3$, then the vertices in $B_i$ are colored one-to-one with the elements in $B_{i,1}$ for each $i\in [3]$, following the order of their indices from smallest to largest.
    \item \label{2.3} If $m_1=m_2$ and $m_2\neq m_3$, then the  vertices in $B_1$, $B_2$ and $B_3$ are colored one-to-one with the elements in $B_{1,1}\cup B_{1,2}$, $B_{2,1}\cup B_{2,2}$ and $B_{3,1}$, respectively, following the order of their indices from smallest to largest.
    \item \label{2.4} If $m_1\neq m_2$ and $m_2=m_3$, then the vertices in $B_1$, $B_2$ and $B_3$ are colored one-to-one with the elements in $B_{1,1}\cup B_{1,4}$, $B_{2,1}$ and $B_{3,1}$, respectively, following the order of their indices from smallest to largest.
\end{enumerate}
\end{comment}

\medskip

\textbf{Case 1: $m_1\leq 1+{3s \over 2} \lceil{\omega\over{\ell-1}}\rceil$ and $\alpha \equiv 0\pmod{3}$}.

For each $i\in [k]$, if $m_i\leq\alpha$,
by Definition \ref{def:balanced_coloring}, we have $L^c_{i, 2(s-1)}=\{1, \cdots, \omega-1\}$ and 
$L^c_{i, 2s-1}=\{\lceil{\ell\omega\over{\ell-1}}\rceil, \cdots, \lceil{\ell\omega\over{\ell-1}}\rceil-\omega+2\}$.
Since $A^c_0=\{1, \cdots, \omega-1\}$, this coloring is proper.

For each $i\in [k]$, if $m_i>\alpha$,
since $m_1\leq 1+{3s \over 2} \lceil{\omega\over{\ell-1}}\rceil$, we have
$L_{i,2(s-1),1}=X_1\cup\{\alpha +i + 3j: 0 \le j \le m_i - \alpha - 1\}$.
If $m_i<{{4-i}\over 2}+{{3s-2}\over 2}\lceil{\omega\over{\ell-1}}\rceil$, then we have $L^c_{i, 2s-1}=\{\lceil{\ell\omega\over{\ell-1}}\rceil, \cdots, \lceil{\ell\omega\over{\ell-1}}\rceil-\omega+2\}$.
Otherwise, by Definition \ref{def:balanced_coloring}, we have
$$L_{i,2s-1,1}=\{\lceil{\ell\omega\over{\ell-1}}\rceil,\cdots, 3m_i-2\alpha+i-2\}\cup (\{3m_i-2\alpha+i-4,\cdots, \beta+g_i\}\backslash \{3j+i: \omega\geq j\geq 0\})$$
where $g_1=2, g_2=g_3=1$ if $\beta\equiv 0\pmod{3}$, and
$g_1=g_2=1, g_3=0$ otherwise.
And we have $L_{i,2s-1,2}=\{\lceil{\ell\omega\over{\ell-1}}\rceil, \cdots,1\}\backslash L_{i,2s-1,1}$.
Since $A^c_0=\{1, \cdots, \omega-1\}$, this coloring is proper.
Therefore, $G$ is $\lceil{\ell\omega\over{\ell-1}}\rceil$-colorable.

\medskip

\textbf{Case 2: $m_1\leq 1+{3s \over 2} \lceil{\omega\over{\ell-1}}\rceil$ and $\alpha \equiv 2\pmod{3}$}.
So $\beta\equiv 0\pmod{3}$.

For each $i\in [k]$, if $m_i\leq\alpha$,
by Definition \ref{def:balanced_coloring}, we have $L^c_{i, 2(s-1)}=X_i\cup (\{\lceil{\omega\over{\ell-1}}\rceil,\cdots,1\}\backslash{X_i})$ and 
$L^c_{i, 2s-1}=\{\lceil{\ell\omega\over{\ell-1}}\rceil, \cdots, \lceil{\ell\omega\over{\ell-1}}\rceil-\omega+2\}$.
Since $A^c_0=\{1, \cdots, \omega-1\}$, this coloring is proper.

For each $i\in [k]$, if $m_i>\alpha$,
since $m_1\leq 1+{3s \over 2} \lceil{\omega\over{\ell-1}}\rceil$, we have
$L_{i,2(s-1),1}=X_i\cup\{\alpha +i +1+ 3j: 0 \le j \le m_i - \alpha - 1\}$.
If $m_i<{{5-i}\over 3}+{{3s-2}\over 2}\lceil{\omega\over{\ell-1}}\rceil$, then we have $L^c_{i, 2s-1}=\{\lceil{\ell\omega\over{\ell-1}}\rceil, \cdots, \lceil{\ell\omega\over{\ell-1}}\rceil-\omega+2\}$.
Otherwise, by Definition \ref{def:balanced_coloring}, we have
$$L_{i,2s-1,1}=\{\lceil{\ell\omega\over{\ell-1}}\rceil,\cdots, 3m_i-2\alpha+i-1\}\cup (\{3m_i-2\alpha+i-3,\cdots, \beta+1\}\backslash \{3j+i: \omega\geq j\geq 0\}).$$

And we have $L_{i,2s-1,2}=\{\lceil{\ell\omega\over{\ell-1}}\rceil, \cdots,1\}\backslash L_{i,2s-1,1}$.
Since $A^c_0=\{1, \cdots, \omega-1\}$, this coloring is proper.
Therefore, $G$ is $\lceil{\ell\omega\over{\ell-1}}\rceil$-colorable.

\medskip

\textbf{Case 3: $m_1>1+{3s \over 2} \lceil{\omega\over{\ell-1}}\rceil$ and $\beta \equiv 0\pmod{3}$}.

Let $g=0$ if $\alpha \equiv 0\pmod{3}$, and $g=1$ otherwise. 
Since $m_1>1+{3s \over 2} \lceil{\omega\over{\ell-1}}\rceil$,
$L_{1,2(s-1),1}=X_1\cup\{\alpha+1+g+3j: 0\leq j\leq {{\beta-\alpha-g-3}\over 3}\}\cup\{\beta+1,\cdots, m_1+\lceil{\omega\over{\ell-1}}\rceil\}$
and 
by Definition \ref{def:balanced_coloring},
$L^c_{1,2s-1}=\{\lceil{\ell\omega\over{\ell-1}}\rceil, \cdots, \lceil{\ell\omega\over{\ell-1}}\rceil-\omega+2\}$.
Since $A^c_0=\{1, \cdots, \omega-1\}$, this coloring is proper.

For each $i\in\{2,3\}$, let $X=X_1$ if $\alpha \equiv 0\pmod{3}$, and $X=X_i$ otherwise.
If $m_i\leq \alpha$, $L^c_{i,2s-1}=\{\lceil{\ell\omega\over{\ell-1}}\rceil, \cdots, \lceil{\ell\omega\over{\ell-1}}\rceil-\omega+2\}$.
If $\alpha<m_i\leq |B'_i|+(s-1)\lceil{\omega\over{l-1}}\rceil$,
$L_{i,2(s-1),1}=X\cup \{\alpha+i+g+3j: 0\leq j\leq m_i-{{\alpha+g+3}\over 3}-(s-1)\lceil{\omega\over{\ell-1}}\rceil\}$
and $L^c_{i,2s-1}=\{\lceil{\ell\omega\over{\ell-1}}\rceil, \cdots, \lceil{\ell\omega\over{\ell-1}}\rceil-\omega+2\}$.
If $m_i>|B'_i|+(s-1)\lceil{\omega\over{l-1}}\rceil$,
$L_{i,2(s-1),1}=X\cup\{\alpha+i+g+3j: 0\leq j\leq {{\beta-\alpha-g-3}\over 3}\}\cup
\{m_1+s\lceil{\omega\over{\ell-1}}\rceil+i-1+2j: 0\leq j\leq m_i-1-{3s-2\over 2}\lceil{\omega\over{\ell-1}}\rceil\}.$
And by Definition \ref{def:balanced_coloring},
we have 
$L_{i, 2s-1,1}=\{\lceil{\ell\omega\over{\ell-1}}\rceil, \cdots,2m_i+m_1+(i-2)-(2s-2)\lceil{\omega\over{\ell-1}}\rceil\}\cup
\{2m_i+m_1-(2s-2)\lceil{\omega\over{\ell-1}}\rceil-(4-i)-2j: 0\leq j\leq m_i-{3s-2\over 2}\lceil{\omega\over{\ell-1}}\rceil-(4-i)\}\cup\{m_1+s\lceil{\omega\over{\ell-1}}\rceil,\cdots, \beta+1\}$ and 
$L_{i,2s-1,2}=\{\lceil{\ell\omega\over{\ell-1}},\cdots,1\rceil\}\backslash L_{i,2s-1,1}$.
Since $A^c_0=\{1, \cdots, \omega-1\}$, this coloring is proper.
Therefore, $G$ is $\lceil{\ell\omega\over{\ell-1}}\rceil$-colorable.

\medskip

\textbf{Case 4: $m_1>1+{3s \over 2} \lceil{\omega\over{\ell-1}}\rceil$ and $\beta \equiv 2\pmod{3}$}. So $\alpha \equiv 0\pmod{3}.$

Since $m_1>1+{3s \over 2} \lceil{\omega\over{\ell-1}}\rceil$,
$L_{1,2(s-1),1}=X_1\cup\{\alpha+1+3j: 0\leq j\leq {{\beta-\alpha-2}\over 3}\}\cup\{\beta+1,\cdots, m_1+\lceil{\omega\over{\ell-1}}\rceil\}$
and 
by Definition \ref{def:balanced_coloring},
$L^c_{1,2s-1}=\{\lceil{\ell\omega\over{\ell-1}}\rceil, \cdots, \lceil{\ell\omega\over{\ell-1}}\rceil-\omega+2\}$.

For each $i\in\{2,3\}$,
if $m_i\leq \alpha$, $L^c_{i,2s-1}=\{\lceil{\ell\omega\over{\ell-1}}\rceil, \cdots, \lceil{\ell\omega\over{\ell-1}}\rceil-\omega+2\}$.
If $\alpha<m_i\leq |B'_i|+(s-1)\lceil{\omega\over{l-1}}\rceil$,
$L_{i,2(s-1),1}=X_1\cup \{\alpha+i+3j: 0\leq j\leq m_i-{{\alpha+3}\over 3}-(s-1)\lceil{\omega\over{\ell-1}}\rceil\}$
and $L^c_{i,2s-1}=\{\lceil{\ell\omega\over{\ell-1}}\rceil, \cdots, \lceil{\ell\omega\over{\ell-1}}\rceil-\omega+2\}$.
If $m_i>|B'_i|+(s-1)\lceil{\omega\over{l-1}}\rceil$,
$L_{i,2(s-1),1}=X_1\cup\{\alpha+i+3j: 0\leq j\leq {{\beta-\alpha-(3i-4)}\over 3}\}\cup
\{m_1+s\lceil{\omega\over{\ell-1}}\rceil+4-i+2j: 0\leq j\leq m_i-{3s-2\over 2}\lceil{\omega\over{\ell-1}}\rceil+{i-4\over 2}\}.$
And by Definition \ref{def:balanced_coloring},
we have 
$L_{i, 2s-1,1}=\{\lceil{\ell\omega\over{\ell-1}}\rceil, \cdots,2m_i+m_1+(i-2)-(2s-2)\lceil{\omega\over{\ell-1}}\rceil\}\cup
\{2m_i+m_1-(2s-2)\lceil{\omega\over{\ell-1}}\rceil-(4-i)-2j: 0\leq j\leq m_i-{3s-2\over 2}\lceil{\omega\over{\ell-1}}\rceil-{3\over 2}\}\cup\{m_1+s\lceil{\omega\over{\ell-1}}\rceil,\cdots, \beta+3-i\}$ and 
$L_{i,2s-1,2}=\{\lceil{\ell\omega\over{\ell-1}},\cdots,1\rceil\}\backslash L_{i,2s-1,1}$.
Since $A^c_0=\{1, \cdots, \omega-1\}$, this coloring is proper.
Therefore, $G$ is $\lceil{\ell\omega\over{\ell-1}}\rceil$-colorable.
\qed

For the sake of convenience, we denote by $A_0$ the intersection of all directed paths in $T$ starting from $a_0$ to $a_i$ for each $i\in [k]$.
\begin{comment}
\begin{figure}[htp]
    \centering
    \includegraphics[width=0.38\linewidth]{5.11.png}
    \caption{the blow-up of $9$-framework with $m=1$}
    \label{fm_1}
\end{figure}
\end{comment}

\begin{lemma}\label{m_1}
    $m\neq1$.
\end{lemma}
%$L_{i, 2(s-1), 1}=X\cup\{\lceil{3(s-1){\lceil {\omega \over {l-1}}\rceil}\over 2}\rceil+i, \cdots, \lceil{3s{\lceil {\omega \over {l-1}}\rceil}\over 2}\rceil-(3-i)+{g\over 2}\}\cup \{m_1+s\lceil{\omega\over{l-1}}\rceil+(i-1)+{g\over 2},\cdots, 2m_i+m_1-(3-i)-(2s-2)\lceil{\omega\over{l-1}}\rceil\}$
%where the second set whose elements form an arithmetic sequence with that first term being $\lceil{3(s-1){\lceil {\omega \over {l-1}}\rceil}\over 2}\rceil+i$ and the common difference being $3$ and the last set whose elements form an arithmetic sequence with that first term being $m_1+s\lceil{\omega\over{l-1}}\rceil+(i-1)+{g\over 2}$ and the common difference being $2$.
%And for each $i\in \{2, 3\}$,
%$L^c_{i, 2s-1}=\{\lceil{l\omega\over{l-1}}\rceil, %\cdots,2m_i+m_1+(i-2)-(2s-2)\lceil{\omega\over{l-1}}\rceil\}\cup
%\{2m_i+m_1-(2s-2)\lceil{\omega\over{l-1}}\rceil-(4-i),\cdots, m_1+s\lceil{\omega\over{l-1}}\rceil+(4-i)+{g\over 2}\}\cup\{m_1+s\lceil{\omega\over{l-1}}\rceil+{g\over 2},\cdots, \lceil{3s{\lceil {\omega \over {l-1}}\rceil}\over 2}\rceil+1+{g\over 2}\}\cup\{2m_i+m_1-(2s-2)\lceil{\omega\over{l-1}}\rceil-(3-i),\cdots, m_1+s\lceil{\omega\over{l-1}}\rceil+(i-1)+{g\over 2}\}\cup \{\lceil{3s{\lceil {\omega \over {l-1}}\rceil}\over 2}\rceil+{g\over 2}, \cdots, 1\}$
%where the common difference in the second set and the forth set is $-2$ and the common difference in the other sets is $-1$.
%By balanced coloring and Lemma \ref{lem1} (1), $G$ is $\lceil{l\omega\over{l-1}}\rceil$-colorable.\qed
\pf Suppose that $m=1$, then $2\leq k-m\leq 3$.
%During the coloring process, we color the vertices in $S$ according to its true size, while all other parts are colored based on a size of $\omega-1$. 
Note that the case when $m=1$ and $k-m=2$ is the same as the case when $m=0$ and $k=3$.
So we assume that $k-m=3$.
Then either $A_1$ is complete to $A_2$ and $B_1$ is complete to $B_3$, $B_4$, 
or $A_1$ is complete to $A_2$, $A_3$ and $B_1$ is complete to $B_4$. 
In all the cases, we may assume that $A_1$ is complete to $A_2$ and $B_1$ is complete to $B_3$, $B_4$.
Let $B^*$ denote the set of internal vertices on the path from $b_3$ to $b_1$ in $S$ (See Figure \ref{fm_1}; A red line means that two sets obey the ordering, while a black line means two parts are complete.).
%Moreover, the cases with $m = 1$ and $k-m=3$ can ultimately be reduced to a single scenario, that is, $A_1$ is complete to $A_2$ and $B_1$ is complete to $B_3$, $B_4$.

\begin{comment}

\begin{figure}[htp]
    \centering
    \includegraphics[width=0.5\linewidth]{m=1,1.jpg}
    \caption{$\ell$-framework with $m=1$}
\end{figure}
\end{comment}
For simplicity, we continue to use the notation in the case when $m=0$ and $k=3$. For example, $m_1,m_2,m_3$ represent the size of the largest, the second largest and the smallest clique among $B_2$, $B_3$, and $B_4$, respectively.

When $A_1$ is complete to $A_2$ and $B_1$ is complete to $B_3$, $B_4$,
since $|B^*|\geq 0$, 
then $B^*$ is complete to $B_2, B_3, B_4$
and $B^*$ with $B_1$ obeys the ordering.
So $|B_2| \leq |B_1|$, we may assume that $|A_2|\geq |A_1|$ (if this case is $\lceil{\ell\omega\over{\ell-1}}\rceil$-colorable,
then the case when $|A_2|\leq |A_1|$ is also $\lceil{\ell\omega\over{\ell-1}}\rceil$-colorable).

When $m_1\leq 1+{3s\over 2}{\lceil{\omega\over{\ell-1}}\rceil}$,
let $G\backslash {\cal P}_1$ be colored as the case when $m=0$ and $k=3$, 
then the color set of $B^*$ is 
$\{m_1+m_2+m_3+1,\cdots,\omega\}$
and $A^c_2$ is a subset of $L^c_{2, 2s-1}$.
From the coloring in the case when $m=0$ and $k-m=3$, we can see that, at most $\lceil{\omega\over{\ell-1}}\rceil$ colors in set $\{\lceil{\ell\omega\over{\ell-1}}\rceil, \cdots, 1\}$ have been shifted backward.
So we can assume that ${A'}^c_1=L^c_{1,2s-1}=\{\lceil{\ell\omega\over{\ell-1}}\rceil-|A_2|,\cdots, \lceil{\ell\omega\over{\ell-1}}\rceil-|A_2|-|A_1|+1\}$.
In fact, $A^c_1$ is the set composed of the elements from ($|A_2|+1$)-th to the ($|A_1|+|A_2|$)-th in $L^c_{2, 2s-1}$ and the difference between ${A'}^c_1$ and $A^c_1$ is at most $\lceil{\omega\over{\ell-1}}\rceil$ elements.
Let $\alpha=\lceil{3(s-1)\over 2} {\lceil {\omega \over {\ell-1}}\rceil} \rceil $ and let
$h=3$ if $\alpha \equiv 0\pmod{3}$, otherwise $h=2$.

Suppose that $|B_2|=m_f$ for each $f\in [3]$.
We have the color $\lceil{\ell\omega\over{\ell-1}}\rceil-|A_2|\not\in B^*$; for otherwise, 
$B^c_1=\{3j+f:0\leq j\leq \omega-m_1-m_2-m_3+m_f-|A_2|-1\}\cup\{\lceil{\ell\omega\over{\ell-1}}\rceil-|A_2|+1,\cdots, \lceil{\ell\omega\over{\ell-1}}\rceil\}.$
\begin{equation*}
\begin{aligned}
 |{A'}^c_1\cap B^c_1|&\leq \lceil {{1\over 3}(3\omega-3m_1-3m_2-3m_3+3m_f-3|A_2|-3+f-(\lceil{\ell\omega\over{\ell-1}}\rceil-|A_2|-|A_1|+1)+1)}\rceil\\
&\leq \omega-m_1-m_2-m_3+m_f-\lfloor{{1\over 3}\lceil{\ell\omega\over{\ell-1}}\rceil}+{{|A_1|+3-f}\over 3}\rfloor
\leq (s-2)\lceil{\omega\over{\ell-1}}\rceil.
\end{aligned}
\end{equation*}
So $|A^c_1\cap B^c_1|\leq (s-1)\lceil{\omega\over{\ell-1}}\rceil$, 
by Lemma \ref{lem11}, $G$ is $\lceil{\ell\omega\over{\ell-1}}\rceil$-colorable.

\medskip
\textbf{Case 1: $m_1\leq 1+{3s \over 2} \lceil{\omega\over{\ell-1}}\rceil$ and $|B_2|=m_1$.}

Now, we assume that $|B_2|=m_1$.
Similarly, if the color $\lceil{\ell\omega\over{\ell-1}}\rceil-|A_2|$ 
 is in $B_{1,3}$, then
$G$ is $\lceil{\ell\omega\over{\ell-1}}\rceil$-colorable. 
 Suppose that color $\lceil{\ell\omega\over{\ell-1}}\rceil-|A_2|$ 
 is in $B_{i,1}$ for each $i\in [3]$.
 Then $B^c_1=\{3j+1: 0\leq j\leq \omega-{2\over 3}\lceil{\ell\omega\over {\ell-1}}\rceil-{|A_2|\over 3}-{i\over 3}\}\cup
\{\lceil{\ell\omega\over {\ell-1}}\rceil-|A_2|+4-i+3j: 0\leq j\leq m_3-{1\over 3}\lceil{\ell\omega\over {\ell-1}}\rceil+{|A_2|\over 3}-{{6-i}\over 3}\}\cup
\{3m_3+1+2j: 0\leq j\leq m_2-m_3-1\}\cup
\{2m_2+m_3+1,\cdots, \lceil{\ell\omega\over {\ell-1}}\rceil\}.$
Note that $B^c_1$ is not an ordered set.
Then
 $|{A'}^c_1\cap B^c_1|\leq \lceil{{1\over 3}(3\omega-2\lceil{\ell\omega\over{\ell-1}}\rceil-|A_2|-i+1-(\lceil{\ell\omega\over{\ell-1}}\rceil-|A_2|-|A_1|+1)+1)}\rceil\leq
 -\lceil{\omega\over{\ell-1}}\rceil+\lceil{{|A_1|+1-i}\over 3}\rceil,$
 and
\begin{flalign*}
 |A^c_1\cap B^c_1|&\leq -\lceil{\omega\over{\ell-1}}\rceil+\lceil{{|A_1|+1-i}\over 3}\rceil+
 \max\{{1\over 3}(3m_1-2\alpha+1-h-\lceil{\ell\omega\over{\ell-1}}\rceil+|A_2|-4+i)+1,0\}\\
& \leq (s-1)\lceil{\omega\over{\ell-1}}\rceil.&
\end{flalign*}
By Lemma \ref{lem11}, $G$ is $\lceil{\ell\omega\over{\ell-1}}\rceil$-colorable.

If the color  $\lceil{\ell\omega\over{\ell-1}}\rceil-|A_2|\in B_{i,2}$ for each $i\in [2]$,
then $B^c_1=\{3j+1:0\leq j\leq \omega-{{1\over 2}\lceil{\ell\omega\over{\ell-1}}\rceil}-{{|A_2|+m_3}\over 2}-{i\over 2}\}\cup\{\lceil{\ell\omega\over{\ell-1}}\rceil-|A_2|+3-i+2j:0\leq j\leq m_2-{{4-i}\over 2}+{{1\over 2}(m_3+|A_2|-\lceil{\ell\omega\over{\ell-1}}\rceil)}\}
\cup \{2m_2+m_3+1,\cdots, \lceil{\ell\omega\over{\ell-1}}\rceil\}.$
Since $\lceil{\ell\omega\over{\ell-1}}\rceil-|A_2|\geq 3m_3+i$ and $m_3\geq s\lceil{\omega\over{\ell-1}}\rceil+1$, then
\begin{equation*}
\begin{aligned}
 |{A'}^c_1\cap B^c_1| &\leq \lceil {{1\over 3}(3\omega-{{3\over 2}\lceil{\ell\omega\over{\ell-1}}\rceil}-{3(|A_2|+m_3)\over 2}-{3i\over 2}+1-(\lceil{\ell\omega\over{\ell-1}}\rceil-|A_2|-|A_1|+1)+1)}\rceil\\
&=\lceil{{1\over 3}({3\omega-{{5\over 2}\lceil{\ell\omega\over{\ell-1}}\rceil}-{|A_2|\over 2}+|A_1|-{3m_3\over 2}-{3i\over 2}+1})}\rceil\\
%\leq 
%\lceil{{3\omega-{5\lceil{l\omega\over{l-1}}\rceil\over 2}+{|A_2|\over 2}-{3m_3\over 2}-{3i\over 2}+1}\over 3}\rceil\\
&\leq \lceil{{1\over 3}({3\omega-{{5\over 2}\lceil{\ell\omega\over{\ell-1}}\rceil}+{{1\over 2}\lceil{l\omega\over{l-1}}\rceil}-{3m_3\over 2}-{i\over 2}-{3m_3\over 2}-{3i\over 2}+1})}\rceil\\
&\leq \omega-m_3-{{2\over 3}\lceil{\ell\omega\over{\ell-1}}\rceil}+{{3-2i}\over 3}
\leq (s-2)\lceil{\omega\over{\ell-1}}\rceil.
\end{aligned}
\end{equation*}
So $|A^c_1\cap B^c_1|\leq (s-1)\lceil{\omega\over{\ell-1}}\rceil$, 
by Lemma \ref{lem11}, $G$ is $\lceil{\ell\omega\over{\ell-1}}\rceil$-colorable.

\medskip
\textbf{Case 2: $m_1\leq 1+{3s \over 2} \lceil{\omega\over{\ell-1}}\rceil$ and $|B_2|=m_2$.}

Now, we assume that $|B_2|=m_2$.
If the color $\lceil{\ell\omega\over{\ell-1}}\rceil-|A_2| \in B_{1,3}$,
then $B^c_1=\{3j+2:0\leq j\leq m_2-{\lceil{\ell\omega\over{\ell-1}}\rceil}-1\}
\cup \{m_1+m_2+m_3+1,\cdots, \lceil{\ell\omega\over{\ell-1}}\rceil\}$
and
\begin{equation*}
\begin{aligned}
 |{A'}^c_1\cap B^c_1| &\leq \lceil {{1\over 3}(-{3\lceil{\omega\over{\ell-1}}\rceil}-1+3m_2-(\lceil{\ell\omega\over{\ell-1}}\rceil-|A_2|-|A_1|+1)+1)}\rceil\\
&\leq -{\lceil{\omega\over{\ell-1}}\rceil}+m_2-{{1\over 3}\lceil{\ell\omega\over{\ell-1}}\rceil}+{{1\over 3}(|A_1|+|A_2|+1)}
\leq (s-2)\lceil{\omega\over{\ell-1}}\rceil.
\end{aligned}
\end{equation*}
So $|A^c_1\cap B^c_1|\leq (s-1)\lceil{\omega\over{\ell-1}}\rceil$, 
by Lemma \ref{lem11}, $G$ is $\lceil{\ell\omega\over{\ell-1}}\rceil$-colorable.

Suppose that the color $\lceil{\ell\omega\over{\ell-1}}\rceil-|A_2|\in B_{i,1}$ for each $i\in [3]$.
 Then $B^c_1=\{3j+2: 0\leq j\leq \omega-{2\over 3}\lceil{\ell\omega\over {\ell-1}}\rceil-{|A_2|\over 3}-{{5-x}\over 3}\}\cup
\{\lceil{\ell\omega\over {\ell-1}}\rceil-|A_2|+x+3j: 0\leq j\leq m_3-{1\over 3}\lceil{\ell\omega\over {\ell-1}}\rceil+{|A_2|\over 3}-{{1+x}\over 3}\}\cup
\{3m_3+2+2j: 0\leq j\leq m_2-m_3-1\}\cup
\{m_1+m_2+m_3+1,\cdots, \lceil{\ell\omega\over {\ell-1}}\rceil\}.$
Then we have
\begin{equation*}
 |{A'}^c_1\cap B^c_1|\leq \lceil{{1\over 3}(3\omega-2\lceil{\ell\omega\over{\ell-1}}\rceil-|A_2|+x-3-(\lceil{\ell\omega\over{\ell-1}}\rceil-|A_2|-|A_1|+1)+1)}\rceil\leq
 -\lceil{\omega\over{\ell-1}}\rceil+\lceil{{|A_1|+x-3}\over 3}\rceil,    
\end{equation*}
while $i=1$, $x=1$ (and respectively $i=2$, $x=3$, and $i=3$, $x=2$) and 
\begin{flalign*}
 |A^c_1\cap B^c_1|&\leq -\lceil{\omega\over{\ell-1}}\rceil+\lceil{{|A_1|+x-3}\over 3}\rceil+\max\{{1\over 3}(3m_2-2\alpha+2-h-\lceil{\ell\omega\over{\ell-1}}\rceil+|A_2|-x)+1,0\}\\
&\leq (s-1)\lceil{\omega\over{\ell-1}}\rceil.&
\end{flalign*}
By Lemma \ref{lem11}, $G$ is $\lceil{\ell\omega\over{\ell-1}}\rceil$-colorable.

If the color $\lceil{\ell\omega\over{\ell-1}}\rceil-|A_2|\in B_{i,2}$ for each $i\in [2]$,
then $B^c_1=\{3j+2:0\leq j\leq \omega-{{1\over 2}\lceil{\ell\omega\over{\ell-1}}\rceil}-{{|A_2|+m_3}\over 2}-{{4-i}\over 2}\}\cup\{\lceil{\ell\omega\over{\ell-1}}\rceil-|A_2|+i+2j:0\leq j\leq m_2-{i\over 2}+{{1\over 2}(m_3+|A_2|-\lceil{\ell\omega\over{\ell-1}}\rceil)}\}
\cup \{m_1+m_2+m_3+1,\cdots, \lceil{\ell\omega\over{\ell-1}}\rceil\}.$
Since $\lceil{\ell\omega\over{\ell-1}}\rceil-|A_2|\geq 3m_3+i$ and $m_3\geq s\lceil{\omega\over{\ell-1}}\rceil+1$, we have
\begin{equation*}
\begin{aligned}
 |{A'}^c_1\cap B^c_1|&\leq \lceil {{1\over 3}(3\omega-{{3\over 2}\lceil{\ell\omega\over{\ell-1}}\rceil}-{3(|A_2|+m_3)\over 2}+{3i\over 2}-4-(\lceil{\ell\omega\over{\ell-1}}\rceil-|A_2|-|A_1|+1)+1)}\rceil\\
&\leq \omega-m_3-{{2\over 3}\lceil{\ell\omega\over{\ell-1}}\rceil}
\leq (s-2)\lceil{\omega\over{\ell-1}}\rceil.
\end{aligned}
\end{equation*}
So $|A^c_1\cap B^c_1|\leq (s-1)\lceil{\omega\over{\ell-1}}\rceil$, 
by Lemma \ref{lem11}, $G$ is $\lceil{\ell\omega\over{\ell-1}}\rceil$-colorable.

\medskip
\textbf{Case 3: $m_1\leq 1+{3s \over 2} \lceil{\omega\over{\ell-1}}\rceil$ and $|B_2|=m_3$.}

Finally, we assume that $|B_2|=m_3$.
If the color $\lceil{\ell\omega\over{\ell-1}}\rceil-|A_2|\in \{3m_3+1,\cdots,m_1+m_2+m_3\}$,
then $B^c_1=\{3j+3:0\leq j\leq m_3-{\lceil{\ell\omega\over{\ell-1}}\rceil}-1\}
\cup \{m_1+m_2+m_3+1,\cdots, \lceil{\ell\omega\over{\ell-1}}\rceil\}$
and
\begin{equation*}
\begin{aligned}
 |{A'}^c_1\cap B^c_1|&\leq \lceil {{1\over 3}(-{3\lceil{\omega\over{\ell-1}}\rceil}+3m_3-(\lceil{\ell\omega\over{\ell-1}}\rceil-|A_2|-|A_1|+1)+1)}\rceil\\
&\leq -{\lceil{\omega\over{\ell-1}}\rceil}+m_3-{{1\over 3}\lceil{\ell\omega\over{\ell-1}}\rceil}+{{1\over 3}(|A_1|+|A_2|+2)}
\leq (s-2)\lceil{\omega\over{\ell-1}}\rceil.
\end{aligned}
\end{equation*}
So $|A^c_1\cap B^c_1|\leq (s-1)\lceil{\omega\over{\ell-1}}\rceil$, 
by Lemma \ref{lem11}, $G$ is $\lceil{\ell\omega\over{\ell-1}}\rceil$-colorable.

Suppose that the color $\lceil{\ell\omega\over{\ell-1}}\rceil-|A_2| \in B_{i,1}$ for each $i\in [3]$.
 Then $B^c_1=\{3j+3: 0\leq j\leq \omega-{2\over 3}\lceil{\ell\omega\over {\ell-1}}\rceil-{|A_2|\over 3}-{{6-x}\over 3}\}\cup
\{\lceil{\ell\omega\over {\ell-1}}\rceil-|A_2|+x+3j: 0\leq j\leq m_3-{1\over 3}\lceil{\ell\omega\over {\ell-1}}\rceil+{|A_2|-x\over 3}\}\cup
\{m_1+m_2+m_3+1,\cdots, \lceil{\ell\omega\over {\ell-1}}\rceil\}.$
Then we have
\begin{equation*}
 |{A'}^c_1\cap B^c_1|\leq \lceil{{1\over 3}(3\omega-2\lceil{\ell\omega\over{\ell-1}}\rceil-|A_2|+x-3-(\lceil{\ell\omega\over{\ell-1}}\rceil-|A_2|-|A_1|+1)+1)}\rceil\leq
 -\lceil{\omega\over{\ell-1}}\rceil+\lceil{{|A_1|+x-3}\over 3}\rceil,    
\end{equation*}
while $i=1$, $x=2$ (respectively $i=2$, $x=1$,  and $i=3$, $x=3$) and 
\begin{flalign*}
 |A^c_1\cap B^c_1|&\leq -\lceil{\omega\over{\ell-1}}\rceil+\lceil{{|A_1|+x-3}\over 3}\rceil+\max\{{1\over 3}(3m_3-2\alpha+3-h-\lceil{\ell\omega\over{\ell-1}}\rceil+|A_2|-x)+1,0\}\\
 &\leq (s-1)\lceil{\omega\over{\ell-1}}\rceil.&
\end{flalign*}
By Lemma \ref{lem11}, $G$ is $\lceil{\ell\omega\over{\ell-1}}\rceil$-colorable.

\medskip
Let $\beta=\lceil{3s\over 2} {\lceil {\omega \over {\ell-1}}\rceil} \rceil $.

\textbf{Case 4: $m_1> 1+{3s \over 2} \lceil{\omega\over{\ell-1}}\rceil$ and $|B_2|\neq m_1$.}

When $m_1>1+{{3s\over 2}\lceil{\omega\over{\ell-1}}\rceil}$ and
$|B_2|\neq m_1$,
let $G\backslash {\cal P}_1$ be colored as the case when $m=0$ and $k=3$, then $A^c_2$ is a subset of $L^c_{2, 2s-1}$.
And 
the first $|A_2|$ elements in $L^c_{2,2s-1}$ are used to color $A_2$ and the next $|A_1|$ elements are used to color $A_1$.
The set of elements from  the $(|A_2|+1)$-th position to the $(|A_1|+|A_2|)$-th position
is denoted as $A^c_1$.
Since $m_1>1+{3s\over 2}{\lceil{\omega\over{\ell-1}}\rceil}$ and $m_3\geq s\lceil{\omega\over{\ell-1}}\rceil+1$,
$m_2-m_3<(s-1)\lceil{\omega\over{\ell-1}}\rceil$.
Hence, $A^c_1\cap B_{2,3}=\emptyset$.
We assume that $|B_2|=m_2$.
Either $L_{2, 2s-1,1}=\{\lceil{\ell\omega\over{\ell-1}}\rceil, \cdots,2m_2+m_1-(2s-2)\lceil{\omega\over{\ell-1}}\rceil\}\cup
\{2m_2+m_1-(2s-2)\lceil{\omega\over{\ell-1}}\rceil-2-2j: 0\leq j\leq m_2-{3s-2\over 2}\lceil{\omega\over{\ell-1}}\rceil-{2+g\over 2}\}\cup\{m_1+s\lceil{\omega\over{\ell-1}}\rceil,\cdots, \beta+1\}$ and 
$L_{2,2s-1,2}=\{\lceil{\ell\omega\over{\ell-1}},\cdots,1\rceil\}\backslash L_{2,2s-1,1}$ where $\beta\equiv 2\pmod{3}$, $g=1$, otherwise 
$g=2$
or $L^c_{2, 2s-1}=\{\lceil{\ell\omega\over{\ell-1}}\rceil, \cdots,\lceil{\omega\over{\ell-1}}\rceil+2\}$.
 Since 
$|A^c_1\cap B^c_2|\leq s\lceil{\omega\over{\ell-1}}\rceil$ and $\{m_1+m_2+m_3+1,\cdots,\lceil{\ell\omega\over{\ell-1}}\rceil\}\cap (A^c_1\cup B^c_2)=\emptyset$,
then $\{m_1+m_2+m_3+1,\cdots,\lceil{\ell\omega\over{\ell-1}}\rceil\}$ can be used to color $B_1$ and
we can select colors from $B^c_2$ to color the remaining vertices in $B_1$ such that $|A^c_1\cap B^c_1|\leq (s-1)\lceil{\omega\over{\ell-1}}\rceil$.
%If $|B_2|=m_3$, it also holds.
By Lemma \ref{lem11}, $G$ is $\lceil{\ell\omega\over{\ell-1}}\rceil$-colorable.

\medskip
\textbf{Case 5: $m_1> 1+{3s \over 2} \lceil{\omega\over{\ell-1}}\rceil$ and $|B_2|=m_1$.}

If $\beta\equiv 2\pmod{3}$, let $g=1$, otherwise, let $g=0$.
And if $\lceil{\ell\omega\over{\ell-1}}\rceil-|A_2|-|A_1|+1>\beta-2+g$, let
$y=s\lceil{\omega\over{\ell-1}}\rceil$.
Otherwise, let
$y=s\lceil{\omega\over{\ell-1}}\rceil-(\lfloor{{\beta-2+g-|A_0|-\lceil{\omega\over{\ell-1}}\rceil-1}\over 3}\rfloor+1)=
s\lceil{\omega\over{\ell-1}}\rceil-\lfloor{{\beta+g-|A_0|-\lceil{\omega\over{\ell-1}}\rceil}\over 3}\rfloor$.

Let $B'_2, B'_3$ and $B'_4$ be colored cyclically where $|B'_2|={{1\over 2}(s\lceil{\omega\over{\ell-1}}\rceil+g)}+y$,  $|B'_3|={{1\over 2}(s\lceil{\omega\over{\ell-1}}\rceil+g)}$ 
and  $|B'_4|={{1\over 2}(s\lceil{\omega\over{\ell-1}}\rceil-g)}$ . Then let $B_3\backslash B'_3$, $B_4\backslash B'_4$ be colored cyclically and finally color $B_2\backslash B'_2$. 

Let $G\backslash {\cal P}_1$ have a balanced coloring.
Then $L^c_{1, 2s-1}=\{\lceil{\ell\omega\over{\ell-1}}\rceil,\cdots,\lceil{\omega\over{\ell-1}}\rceil+2\}$ and $A^c_2=\{\lceil{\ell\omega\over{\ell-1}}\rceil,\cdots, \lceil{\ell\omega\over{\ell-1}}\rceil-|A_2|+1\}$,
$A^c_1=\{\lceil{\ell\omega\over{\ell-1}}\rceil-|A_2|,\cdots, \lceil{\ell\omega\over{\ell-1}}\rceil-|A_2|-|A_1|+1\}$.
Since 
$|A_2|\geq|A_1|$, we have $A^c_1\cap (B_2\backslash B'_2)^c=\emptyset$.
%$|B_2\backslash B'_2|+\lceil{\omega\over{l-1}}\rceil=m_1-\lceil{3s\lceil{\omega\over{l-1}}\rceil\over2}\rceil+\lfloor{{\lceil{3s\lceil{\omega\over{l-1}}\rceil\over2}\rceil+g-|A_0|-\lceil{\omega\over{l-1}}\rceil}\over 3}\rfloor+\lceil{\omega\over{l-1}}\rceil<s\lceil{\omega\over{l-1}}\rceil$,
%the color $\lceil{l\omega\over{l-1}}\rceil-|A_2|$ is used to color $B'_2$ or $B_3$ or $B_4$.
Since $\{m_1+m_2+m_3+1,\cdots,\lceil{\ell\omega\over{\ell-1}}\rceil\}$ can be used to color $B_1$,
%and $\{\lceil{l\omega\over{l-1}}\rceil,\cdots,\omega+1\}\cap (A^c_1\cup B^c_2)=\emptyset$,
we can select some colors from $B^c_2$ to color the remaining vertices in $B_1$ such that $|A^c_1\cap B^c_1|\leq (s-1)\lceil{\omega\over{\ell-1}}\rceil$.
By Lemma \ref{lem11}, $G$ is $\lceil{\ell\omega\over{\ell-1}}\rceil$-colorable.
\qed

\begin{lemma}
    $m\neq2$.
\end{lemma}

\pf If $m=2$, then $k=5$ and
either $A_4$ is complete to $A_1$ or $A_4$ is anticomplete to $A_1$.
When $A_4$ is complete to $A_1$, since $G$ is a minimal counterexample of Lemma \ref{lem-color-framework}, by Lemma \ref{Claim 3.1.1}, $|A_0|\geq s{\lceil {\omega \over {\ell-1}}\rceil}+1$.
And by Definition \ref{def:l-framework}, $\cup^3_{i=0}A_i$ is a clique of size greater than $\omega$, a contradiction.

\begin{comment}

\begin{figure}[htp]
    \centering
    \includegraphics[width=0.5\linewidth]{m=2.jpg}
    \caption{$\ell$-framework with $m=2$}
\end{figure}
\end{comment}
When $A_4$ is anticomplete to $A_1$, we may assume that $|A_3|\geq \max\{|A_1|, |A_2|\}$, then let $G\backslash ({\cal P}_1\cup {\cal P}_2)$ be colored as the case when $m=0$ and $k=3$.
Similarly to the proof of Lemma \ref{m_1}, 
we can obtain that for each $i\in [2]$,
$|A^c_i\cap B^c_i|\leq (s-1)\lceil{\omega\over{\ell-1}}\rceil$.
Then by Lemma \ref{lem11}, $G$ is $\lceil{\ell\omega\over{\ell-1}}\rceil$-colorable.
 \qed

\medskip

Therefore, by the above lemmas, when $\ell \equiv 1 \pmod 4$ and $\ell\ge 9$, the blow-up of $\ell$-frameworks $G$ is $\lceil{\ell\over{\ell-1}}\omega(G)\rceil$-colorable.
Hence, Theorem \ref{main theorem} follows from Lemmas \ref{lem:blow-up-cycle-colorable} and \ref{lem-color-framework}.

\medskip

Finally, we propose the following conjecture.
\begin{conjecture}
If $G$ is a $5$-holed graph, then $\chi(G)\leq \lceil {5 \over 4} \omega(G) \rceil$.
\end{conjecture}

\section{Appendix}
\subsection{Proof of (\ref{3.1})}
In this section, we prove (\ref{3.1}):
\begin{equation*}
\sum^{i'}_{j=1}(m_j-\lceil{sk\lceil{\omega\over{\ell-1}}\rceil\over {k-1}}\rceil)+k\lceil{s\lceil{\omega\over{\ell-1}}\rceil\over {k-1}}\rceil+\omega-m_1\leq \lceil{\ell\omega\over{\ell-1}}\rceil.
\end{equation*}

When $k=5$ and $i'=5$, then $l=7$, and 
\begin{equation*}
\begin{aligned}
&\sum^{i'}_{j=1}(m_j-\lceil{sk\lceil{\omega\over{\ell-1}}\rceil\over {k-1}}\rceil)+k\lceil{s\lceil{\omega\over{\ell-1}}\rceil\over {k-1}}\rceil+\omega-m_1-\lceil{\ell\omega\over{\ell-1}}\rceil\\
&\leq m_1+m_2+m_3+m_4+m_5-{4\times {5(\ell-3)\lceil{\omega\over{\ell-1}}\rceil\over 4}\over 4}+\omega-m_1-\lceil{\ell\omega\over{\ell-1}}\rceil\\
&\leq m_2+m_3+m_4+m_5-{{5(\ell-3)\omega+4\omega}\over 4(\ell-1)}\leq {4\omega\over 5}-{{(5\ell-11)\omega}\over {4(\ell-1)}}<0.
\end{aligned}
\end{equation*}

When $k=5$ and $i'=4$, then 
\begin{equation*}
\begin{aligned}
&\sum^{i'}_{j=1}(m_j-\lceil{sk\lceil{\omega\over{\ell-1}}\rceil\over {k-1}}\rceil)+k\lceil{s\lceil{\omega\over{\ell-1}}\rceil\over {k-1}}\rceil+\omega-m_1-\lceil{\ell\omega\over{\ell-1}}\rceil\\
&\leq m_2+m_3+m_4-{3\times {5(\ell-3){\omega\over{\ell-1}}\over 4}\over 4}-{\omega\over{\ell-1}}
\leq m_2+m_3+m_4-{{(15\ell-29)\omega}\over 16(\ell-1)}\leq {(-3\ell+13)\omega\over 8(\ell-1)}<0.
\end{aligned}
\end{equation*}

When $k=5$ and $i'=3$, then 
\begin{equation*}
\begin{aligned}
&\sum^{i'}_{j=1}(m_j-\lceil{sk\lceil{\omega\over{\ell-1}}\rceil\over {k-1}}\rceil)+k\lceil{s\lceil{\omega\over{\ell-1}}\rceil\over {k-1}}\rceil+\omega-m_1-\lceil{\ell\omega\over{\ell-1}}\rceil\\
&\leq m_2+m_3-{2\times {5(\ell-3){\omega\over{\ell-1}}\over 4}\over 4}-{\omega\over{\ell-1}}
\leq {2(\omega-{\omega\over{\ell-1}}{{\ell-3}\over 4}\times 2)\over 3}-{{(5\ell-15+8)\omega}\over 8(\ell-1)}\leq {(-7\ell+29)\omega\over 24(\ell-1)}<0.
\end{aligned}
\end{equation*}

When $k=5$ and $i'=2$, then 
\begin{equation*}
\begin{aligned}
&\sum^{i'}_{j=1}(m_j-\lceil{sk\lceil{\omega\over{\ell-1}}\rceil\over {k-1}}\rceil)+k\lceil{s\lceil{\omega\over{\ell-1}}\rceil\over {k-1}}\rceil+\omega-m_1-\lceil{\ell\omega\over{\ell-1}}\rceil\\
&\leq m_2-{{5(\ell-3){\omega\over{\ell-1}}\over 4}\over 4}-{\omega\over{\ell-1}}+1
\leq {\omega-{\omega\over{\ell-1}}{{\ell-3}\over 4}\times 3-3\over 2}-{{(5\ell-15+16)\omega}\over 16(\ell-1)}+1\leq {(-3\ell+9)\omega\over 16(\ell-1)}<0.
\end{aligned}
\end{equation*}

When $k=5$ and $i'=1$, then 
\begin{equation*}
\begin{aligned}
&\sum^{i'}_{j=1}(m_j-\lceil{sk\lceil{\omega\over{\ell-1}}\rceil\over {k-1}}\rceil)+k\lceil{s\lceil{\omega\over{\ell-1}}\rceil\over {k-1}}\rceil+\omega-m_1-\lceil{\ell\omega\over{\ell-1}}\rceil\\
&=m_1-\lceil{5\lceil{\omega\over{\ell-1}}\rceil\over 4}\rceil+5\lceil{\lceil{\omega\over{\ell-1}}\rceil\over 4}\rceil+\omega-m_1-\lceil{\ell\omega\over{\ell-1}}\rceil=4\lceil{\lceil{\omega\over{\ell-1}}\rceil\over 4}\rceil-2\lceil{\omega\over{\ell-1}}\rceil.
\end{aligned}
\end{equation*}
If $\lceil{\omega\over{\ell-1}}\rceil=1$, then $\omega\leq \ell-1=6$ which contradicts to $\omega\geq \sum^5_{i=1}m_i\geq 2\times 5=10$.
So $\lceil{\omega\over{\ell-1}}\rceil\geq 2$, and 
$$\sum^{i'}_{j=1}(m_j-\lceil{sk\lceil{\omega\over{\ell-1}}\rceil\over {k-1}}\rceil)+k\lceil{s\lceil{\omega\over{\ell-1}}\rceil\over {k-1}}\rceil+\omega-m_1-\lceil{\ell\omega\over{\ell-1}}\rceil\leq 0.$$

When $k=4$ and $i'=4$, then 
\begin{equation*}
\begin{aligned}
&\sum^{i'}_{j=1}(m_j-\lceil{sk\lceil{\omega\over{\ell-1}}\rceil\over {k-1}}\rceil)+k\lceil{s\lceil{\omega\over{\ell-1}}\rceil\over {k-1}}\rceil+\omega-m_1-\lceil{\ell\omega\over{\ell-1}}\rceil\\
&\leq m_2+m_3+m_4-{3\times {4(\ell-3){\omega\over{\ell-1}}\over 4}\over 3}-{\omega\over{\ell-1}}
\leq {3\omega\over 4}-{{(\ell-2)\omega}\over \ell-1}\leq {(-\ell+5)\omega\over 4(\ell-1)}<0.
\end{aligned}
\end{equation*}

When $k=4$ and $i'=3$, then 
\begin{equation*}
\begin{aligned}
&\sum^{i'}_{j=1}(m_j-\lceil{sk\lceil{\omega\over{\ell-1}}\rceil\over {k-1}}\rceil)+k\lceil{s\lceil{\omega\over{\ell-1}}\rceil\over {k-1}}\rceil+\omega-m_1-\lceil{\ell\omega\over{\ell-1}}\rceil\\
&\leq m_2+m_3-{2\times {4(\ell-3){\omega\over{\ell-1}}\over 4}\over 3}-{\omega\over{\ell-1}}
\leq {{2\omega-{2\omega\over {\ell-1}}{\ell-3\over 4}}\over 3}-{{(2\ell-3)\omega}\over 3(\ell-1)}\leq {(-\ell+5)\omega\over 6(\ell-1)}<0.
\end{aligned}
\end{equation*}

When $k=4$ and $i'=2$, then 
\begin{equation*}
\begin{aligned}
&\sum^{i'}_{j=1}(m_j-\lceil{sk\lceil{\omega\over{\ell-1}}\rceil\over {k-1}}\rceil)+k\lceil{s\lceil{\omega\over{\ell-1}}\rceil\over {k-1}}\rceil+\omega-m_1-\lceil{\ell\omega\over{\ell-1}}\rceil\\
&\leq m_2-{2{4(\ell-3){\omega\over{\ell-1}}\over 4}\over 3}-{\omega\over{\ell-1}}
\leq {{\omega-{2\omega\over {\ell-1}}{\ell-3\over 4}}\over 2}-{{\ell\omega}\over 3(\ell-1)}\leq {(-\ell+3)\omega\over 12(\ell-1)}<0.
\end{aligned}
\end{equation*}

When $k=4$ and $i'=1$, then 
\begin{equation*}
\begin{aligned}
\sum^{i'}_{j=1}(m_j-\lceil{sk\lceil{\omega\over{\ell-1}}\rceil\over {k-1}}\rceil)+k\lceil{s\lceil{\omega\over{\ell-1}}\rceil\over {k-1}}\rceil+\omega-m_1-\lceil{\ell\omega\over{\ell-1}}\rceil
=3\lceil{s\lceil{\omega\over{\ell-1}}\rceil\over 3}\rceil-(s+1)\lceil{\omega\over{\ell-1}}\rceil.
\end{aligned}
\end{equation*}
If $\lceil{\omega\over{\ell-1}}\rceil=1$, then $\omega\leq \ell-1$ which contradicts to $\omega\geq \sum^4_{i=1}m_i\geq (s+1)\times 4=\ell+1$.
So $\lceil{\omega\over{\ell-1}}\rceil\geq 2$, and 
$$\sum^{i'}_{j=1}(m_j-\lceil{sk\lceil{\omega\over{\ell-1}}\rceil\over {k-1}}\rceil)+k\lceil{s\lceil{\omega\over{\ell-1}}\rceil\over {k-1}}\rceil+\omega-m_1-\lceil{\ell\omega\over{\ell-1}}\rceil\leq 0.$$

When $k=3$ and $i'=3$, then 
\begin{equation*}
\begin{aligned}
&\sum^{i'}_{j=1}(m_j-\lceil{sk\lceil{\omega\over{\ell-1}}\rceil\over {k-1}}\rceil)+k\lceil{s\lceil{\omega\over{\ell-1}}\rceil\over {k-1}}\rceil+\omega-m_1-\lceil{\ell\omega\over{\ell-1}}\rceil\\
&\leq m_2+m_3-{2\times {3(\ell-3){\omega\over{\ell-1}}\over 4}\over 2}-{\omega\over{\ell-1}}
\leq {2\omega\over 3}-{{(3\ell-5)\omega}\over 4(\ell-1)}\leq {(-\ell+7)\omega\over 12(\ell-1)}\leq 0.
\end{aligned}
\end{equation*}

When $k=3$ and $i'=2$, then 
\begin{equation*}
\begin{aligned}
&\sum^{i'}_{j=1}(m_j-\lceil{sk\lceil{\omega\over{\ell-1}}\rceil\over {k-1}}\rceil)+k\lceil{s\lceil{\omega\over{\ell-1}}\rceil\over {k-1}}\rceil+\omega-m_1-\lceil{\ell\omega\over{\ell-1}}\rceil\\
&\leq m_2-{ {3(\ell-3){\omega\over{\ell-1}}\over 4}\over 2}-{\omega\over{\ell-1}}
\leq {{\omega-{\omega\over {\ell-1}}{\ell-3\over 4}}\over 2}-{{(3\ell-1)\omega}\over 8(\ell-1)}=<0.
\end{aligned}
\end{equation*}

When $k=3$ and $i'=1$, then 
\begin{equation*}
\begin{aligned}
\sum^{i'}_{j=1}(m_j-\lceil{sk\lceil{\omega\over{\ell-1}}\rceil\over {k-1}}\rceil)+k\lceil{s\lceil{\omega\over{\ell-1}}\rceil\over {k-1}}\rceil+\omega-m_1-\lceil{\ell\omega\over{\ell-1}}\rceil
=2\lceil{s\lceil{\omega\over{\ell-1}}\rceil\over 2}\rceil-(s+1)\lceil{\omega\over{\ell-1}}\rceil\leq 0.
\end{aligned}
\end{equation*}

\end{document}